\documentclass[leqno,12pt]{article}
\usepackage{amsmath,amssymb}
\usepackage{color}
\usepackage{epsfig}
\usepackage[utf8]{inputenc}
\usepackage[T1]{fontenc}

\newtheorem{Theorem}{\hspace{\parindent}\bf Theorem}[section]
\newtheorem{Lemma}{\hspace{\parindent}\bf Lemma}[section]
\newtheorem{Proposition}{\hspace{\parindent}\bf Proposition}[section]
\newtheorem{Corollary}{\hspace{\parindent}\bf Corollary}[section]

\newtheorem{Definition}{\hspace{\parindent}\bf Definition}[section]

\renewcommand{\theequation}{\arabic{section}.\arabic{equation}}

\newcommand{\qed}{\hfill$\square$\vspace{0.3cm}}

\newcommand{\R}{\mathbb{R}}
\newcommand{\ren}{\mathbb{R}^N}
\newcommand{\Tr}{\mathop{\rm Tr}}
\newcommand{\E}{\mathop{\rm E}}

\begin{document}

\title{\textbf{Classical solutions for a logarithmic fractional diffusion equation}}
\author{ by \\
Arturo de Pablo, Fernando Quir\'{o}s, \\ Ana Rodr\'{\i}guez, and Juan Luis V\'{a}zquez}

\maketitle

\begin{abstract}
We prove global existence and uniqueness of  strong solutions to the
logarithmic porous medium type equation with fractional diffusion
$$
\partial_tu+(-\Delta)^{1/2}\log(1+u)=0,
$$
posed for $x\in \mathbb{R}$, with nonnegative initial data in some
function space of $L \log\!L$ type. The solutions are shown to
become bounded and $C^\infty$  smooth in $(x,t)$ for all positive times.
We also reformulate this equation as a transport equation with nonlocal velocity and
critical viscosity,   a topic of current relevance. Interesting functional inequalities are involved.
\end{abstract}

\vskip 1cm

\noindent{\makebox[1in]\hrulefill}\newline

%

\section{Introduction}\label{sect-introduction}
\setcounter{equation}{0}

In this paper we develop the basic existence, uniqueness and
regularity theory for  the problem
\begin{equation}  \label{eq:main}
\left\{
\begin{array}{ll}
\partial_tu+(-\Delta)^{1/2}\log(1+u)=0, & \qquad  x\in\mathbb{R},\; t>0,
\\ [4mm]
u(x,0) = f(x)\ge0, & \qquad x\in\mathbb{R}.%
\end{array}
\right.
\end{equation}
The equation in \eqref{eq:main} can be viewed as the limit $m\to0$
 in the so-called fractional porous  medium equation,
\begin{equation}
\label{eq:fpme}
\partial_tu+(-\Delta)^{\sigma/2}u^m=0, \qquad m>0, \quad 0<\sigma<2,
\end{equation}
after a shift in the $u$-variable and a change in the time scale.
The latter equation was treated in our papers \cite{pqrv},
\cite{pqrv2}, where it was proved that it generates a contraction
semigroup in $L^1(\ren)$  for any dimension $N\ge1$, and that
solutions become instantaneously bounded and $C^\alpha$ in space and
time for data in $L^1(\ren)\cap L^p(\ren)$ with $p\ge 1$ larger than
a critical value $p_*=N(1-m)/\sigma$.

The difficulty we face here is that, according to those papers, the
logarithmic diffusion is borderline for regularity questions when
$\sigma=1$ and $N=1$ for  data in the natural space $L^1(\R)$. This
entails a very delicate critical-case analysis and a new type of
regularity results: besides the expected result for $f\in
L^1(\R)\cap L^p(\R)$ with $p>1$, we obtain that solutions become
immediately bounded when $f$ belongs to an $L \log\!L$ space, almost
$L^1(\mathbb{R})$ but not quite. This also offers some novelty when
compared to the existing results for the standard  porous medium
equation, given by~\eqref{eq:fpme} with $\sigma=2$, which are
gathered in \cite{JLVSmoothing} and \cite{vazquez}. Actually, we go
on to prove that the solutions are  $C^\infty$ in space and time, hence
classical.

Let us remark that the method proposed to tackle regularity has a
more general scope. Actually, it can be applied to positive solutions of
equations of the form $\partial_t u+(-\Delta)^{\sigma/2}\varphi (u)=0$
posed in $\mathbb{R}^N$ under quite unrestrictive assumptions on the
nonlinearity. We will study this issue in a forthcoming work.

A further motivation for our study comes from the following connection:
equation \eqref{eq:main} can be transformed through a special
nonlocal change of variables of the B\"acklund type  into the
transport equation
\begin{equation}
\partial_\tau v-\widetilde H(v)\partial_y v+\partial_y\widetilde H(v)=0,\qquad y\in\mathbb{R},\quad \tau>0,
\label{eq:transport-rare}
\end{equation}
 where $\widetilde H$ stands
for a nonlocal operator which is a modification of the Hilbert
transform. If instead of $\widetilde H$ we had the standard Hilbert
transform $H$, using the identity $\partial_y H=(-\Delta)^{1/2}$,
valid when these operators are applied to regular  functions, we
would get
 \begin{equation} \label{eq:transport}
 \partial_\tau u -H(v)\,\partial_yv + (-\Delta)^{1/2}v=0,
 \end{equation}
which is the transport equation with fractional diffusivity proposed
by C\'ordoba, C\'{o}rdoba and Fontelos in~\cite{ccf}. This is one of
the several one-dimensional models considered in the last years to
recast the main properties of the three-dimensional incompressible
Euler equation and the two-dimensional quasigeostrophic equation,
beginning with the work by Constantin, Lax and
Majda~\cite{Constantin-Lax-Majda}.

Conveniently reformulated, our results for problem~\eqref{eq:main}
produce existence and uniqueness of a \emph{classical} global in
time solution for equation \eqref{eq:transport-rare} for all initial
data in $L^1(\mathbb{R})$. This is a remarkable variation with
respect to the results available for problem~\eqref{eq:transport}: a
global in time solution is only known to exist if the initial value
belongs to the Sobolev space of fractional order
$H^{1/2}(\mathbb{R})$, in which case it is in $H^1(\mathbb{R})$ for
almost every $t>0$; see  Dong~\cite{Dong-2008}. From this regularity one
might try to use the techniques of Kiselev, Nazarov and Shterenberg~\cite{Kiselev-Nazarov-Shterenberg} to obtain further smoothness in
space.

The application of the present approach to the transport
equation~\eqref{eq:transport} is not immediate  and needs further
study. We believe that the connection between fractional diffusion
and nonlocal transport problems is worth pursuing, since it may lead
to the fruitful combination of very different techniques.

As  said above,  the case $\sigma=N=1$ is critical in various
aspects, in particular with respect to Sobolev embeddings. Thus, in
the course of the proof of the smoothing effect in the mentioned
$L\log\!L$ space we need to use a critical fractional Trudinger type
embedding due to Strichartz; see \cite{Strichartz-1972}. We
generalize this embedding to other values of the exponents, a result
in pure functional analysis that we hope could be of further
application.

\section{Preliminaries and main results}\label{sect-main.results}
\setcounter{equation}{0}

We recall that the nonlocal operator $(-\Delta)^{\sigma/2}$,
$\sigma\in(0,2)$,  is defined for any function
$g:\mathbb{R}^N\to\mathbb{R}^N$ in the Schwartz class through the
Fourier transform,
\begin{equation*}
\label{def-fourier}
\mathcal{F}\left((-\Delta)^{\sigma/2}  g\right)(\xi)=|\xi|^{\sigma}
\mathcal{F}(g)(\xi),
\end{equation*}
or via the (hypersingular) Riesz potential,
\begin{equation}
\label{def-riesz}
(-\Delta)^{\sigma/2}  g(x)= C_{N,\sigma}\mbox{
P.V.}\int_{\mathbb{R}} \frac{g(x)-g(y)}{|x-y|^{N+\sigma }}\,dy,
\end{equation}
where $C_{N,\sigma}$ is a normalization constant; see for example
\cite{Landkof}. In our case, $N=\sigma=1$, the constant is
$C_{1,1}=1/\pi$, and we also have $(-\Delta)^{1/2}=H\partial_x $,
where $H$ denotes the Hilbert transform operator, defined trough
\begin{equation*} Hf(x)=\frac1\pi\mbox{
P.V.}\int_{\mathbb{R}}\frac{f(y)}{x-y}\,dy.
\end{equation*}

If we multiply  the equation in \eqref{eq:main} by a test function
$\varphi$ and \lq\lq integrate by parts'',  we  obtain
\begin{equation}\label{weak-nonlocal}
\displaystyle \int_0^\infty\int_{\mathbb{R}}u\,\partial_t
\varphi\,dxdt-\int_0^\infty\int_{\mathbb{R}}(-\Delta)^{1/4}\log(1+u)(-\Delta)^{1/4}\varphi\,d
xdt=0.
\end{equation}
This identity will be the  basis of our definition of a weak
solution. The integrals in~\eqref{weak-nonlocal} make sense if $u$
and $\log(1+u)$ belong to suitable spaces. The  correct space  for
 $\log(1+u)(\cdot,t)$ is the homogeneous fractional Sobolev space
$\dot{H}^{1/2}(\mathbb{R})$, defined as the completion of
$C_0^\infty(\mathbb{R})$ with the norm
$$
  \|\psi\|_{\dot{H}^{1/2}}=\left(\int_{\mathbb{R}}
|\xi||\hat{\psi}|^2\,d\xi\right)^{1/2} =\|(-\Delta)^{1/4}\psi\|_{2}.
$$
The Sobolev space $H^{1/2}(\mathbb{R})$ is then defined through the
norm
$$
  \|\psi\|_{H^{1/2}}=\|\psi\|_{2}+\|(-\Delta)^{1/4}\psi\|_{2}.
$$

\begin{Definition}\label{def:weak.solution.nonlocal} A function $u$
is a {\it weak} $L^1$-energy solution to problem \eqref{eq:main} if:
\begin{itemize}
\item $u\in C([0,\infty): L^1(\mathbb{R}))$ \ and \ $\log(1+u)\in L^2((0,T):\dot{H}^{1/2}(\mathbb{R}))$ for every $T>0$;
\item  identity  \eqref{weak-nonlocal}
holds for every $\varphi\in C_0^1(\mathbb{R}\times(0,\infty))$;
\item  $u(\cdot,0)=f$ almost everywhere.
\end{itemize}
\end{Definition}
For the sake of brevity, we will denote the solutions obtained below
according to this definition merely as weak solutions. We remark
that this is not the only way of defining a solution to
problem~\eqref{eq:main}. There are other possibilities, for instance
entropy solutions \cite{CJ}, useful when dealing with equations
involving convection terms.

As for the initial data, our concept of solution only requires in
principle $f\in L^1_+(\mathbb{R})$. However, in order to prove
existence we will ask  $f$  to belong to the slightly smaller
$L\log\!L$-type space
\begin{equation*}
{\cal X}=\left\{f\ge0  \mbox{
measurable}:\,\int_{\mathbb{R}}(1+f)\log(1+f)<\infty\right\}.
\end{equation*}
Notice that $L_+^1(\mathbb{R})\cap L^p(\mathbb{R})\subset{\cal
X}\subset L_+^1(\mathbb{R})$ for any $p>1$. This is to be compared
with the result for the fractional porous medium equation
\eqref{eq:fpme}, where in the critical case $\sigma=N(1-m)$ we have
required $f\in L^1(\mathbb{R}^N)\cap L^p(\mathbb{R}^N)$ for some
$p>1$.

The space ${\cal X}$ is  natural for problem \eqref{eq:main}.
Indeed, let $\Psi(s)=(1+s)\log(1+s)-s$. A function $f$ belongs to
$\mathcal{X}$ if and only if $f\in L_+^1(\mathbb{R}^N)$ and
$\int_{\mathbb{R}}\Psi(f)<\infty$. On the other hand, after an
integration by parts we formally obtain
$$
\int_{\mathbb{R}}\Psi(u(\cdot,t))=\int_{\mathbb{R}}\Psi(f)
-\int_0^t\int_{\mathbb{R}}\left|(-\Delta)^{1/4}(\log(1+u))\right|^2\le
\int_{\mathbb{R}}\Psi(f),
$$
and we conclude that the space $\mathcal{X}$ is preserved by the
evolution.

\noindent\emph{Notation. } We will denote $
L_{\mathcal{X}}(f):=\int_{\mathbb{R}}\Psi(f)=\int_{\mathbb{R}}\left((1+f)\log(1+f)-f\right)$.

\medskip

Though we will be able to prove existence of a weak solution for any
$f\in\mathcal{X}$, in order to prove uniqueness we will
restrict ourselves to the smaller class of \emph{strong solutions}.

\begin{Definition} We say that a weak solution $u$  to
problem \eqref{eq:main} is a \emph{strong solution} if \
$\partial_tu\in L^1_{{\rm loc}}((0,\infty)\times\mathbb{R})$.
\label{def:strong.solution}
\end{Definition}
If $u$ is a strong solution, then $(-\Delta)^{1/2}\log(1+u)$ is also
an $L^1_{\textrm{loc}}$-function and the equation in~\eqref{eq:main}
is satisfied a.e.

Our first result shows that  problem~\eqref{eq:main} is well posed
in the class of strong solutions  for initial data in $\mathcal{X}$.

\begin{Theorem}\label{th:existence} For every $f\in
{\cal X}$ there exists a unique  strong  solution to problem~\eqref{eq:main}.
\end{Theorem}

Existence and uniqueness use an alternative formulation of
problem~\eqref{eq:main} based in the Dirichlet to Neumann operator.
Given a smooth bounded function $g:\mathbb{R}\mapsto\mathbb{R}$, we
define its harmonic extension
 $v=\E(g)$ to the upper half-plane $\mathbb{R}^2_+$ as the unique smooth bounded solution to
\begin{equation*}
\left\{
\begin{array}{ll}
\Delta_{x,y}  v=0,\qquad &x\in\mathbb{R},\, y>0,\\
v(x,0)=g(x),\qquad&x\in\mathbb{R}.
\end{array}
\right. \label{alfa-extension}
\end{equation*}
Then, it turns out that $-\partial_y v(x,0)=(-\Delta_x)^{1/2} g(x)$,
where $\Delta_{x,y}$ is the Laplacian in all $(x,y)$-variables and
$\Delta_x$ acts only on the $x$-variables (in the sequel we will
drop the subscripts when no confusion arises). The extension
operator $\E$ can be defined by density in the space $\dot
H^{1/2}(\mathbb{R})$, and it is an isometry between this space and
the space  $\mathcal{H}$ defined as the completion of
$C_0^\infty(\overline{\mathbb{R}^2_+})$ with the norm
$$
\|\psi\|_{\mathcal{H}}=\left(\int_0^\infty\int_{\mathbb{R}}|\nabla
\psi|^2\right)^{1/2}.
$$
Therefore,
\begin{equation}
\label{eq:plancherel}
\int_{\mathbb{R}}(-\Delta)^{1/4}\phi\,(-\Delta)^{1/4}\psi=
\int_0^\infty\int_{\mathbb{R}}\nabla \E(\phi)\cdot\nabla\E (\psi).
\end{equation}
We also have
\begin{equation}
\label{eq:fundamental.property.extensions}
\int_0^\infty\int_{\mathbb{R}}\nabla \E(\phi)\cdot\nabla\E (\psi)=
\int_0^\infty\int_{\mathbb{R}}\nabla \eta\cdot\nabla\E (\psi).
\end{equation}
for any $\eta\in\mathcal{H}$ such that $\Tr(\eta)=\phi$; see
\cite{pqrv2}.

Using this approach, problem~\eqref{eq:main} can be written in an
equivalent local form. If $u$ is a solution, then  $w=\E(\log(1+u))$
solves
\begin{equation}
\label{pp:local}
\left\{
\begin{array}{ll}
\Delta w=0,\qquad &(x,y)\in\mathbb{R}^{2}_+,\, t>0,\\
\partial_y w-\partial_t\beta(w)=0,\qquad&x\in\mathbb{R},\,y=0,\, t>0,\\
w=\log(1+f),\qquad&x\in\mathbb{R}, \,y=0,\, t=0,
\end{array}
\right. \qquad\beta(w)=e^{w}-1.
\end{equation}
Conversely, if we obtain a solution $w$
to~\eqref{pp:local}, then $u=\beta(w)\big|_{y=0}$ is a solution
to~\eqref{eq:main}.

We next state the main properties of the solution obtained in the
paper.

\begin{Theorem}\label{th:props} Let $f\in
{\cal X}$. The unique  strong  solution $u$ to problem
\eqref{eq:main} satisfies:
\begin{enumerate}
\item[\rm (i)]   $ \partial_tu\in
L^2(\mathbb{R}\times(\tau,\infty))$ for all $\tau>0$;
\item[\rm (ii)] \emph{$\mathcal{X}$--$L^\infty$ smoothing effect:}
\begin{equation}\label{eq:full.smoothing.effect}
\|u(\cdot,t)\|_\infty\le C\max\{
t^{-1}\mbox{\rm exp}(Ct^{-1/2}\big(L_{\mathcal{X}}(f)\big)^{1/2}),
t^{-3/4}\|f\|_1^{1/2}\big(L_{\mathcal{X}}(f)\big)^{1/4}\} ;
\end{equation}
\item[\rm (iii)] $L_{\mathcal{X}}(u(\cdot,t))$ and  $\|u(\cdot,t)\|_p$, $1\le
p\le\infty$, are non-increasing functions of $t$ in $(0,\infty)$;
\item[\rm (iv)] $\displaystyle\int_{\mathbb{R}}u(x,t)\,dx=
\int_{\mathbb{R}}f(x)\,dx$ for every $t\ge0$ \emph{(conservation of
mass)};
\item[\rm (v)]  $u\in C^\infty(\mathbb{R}\times(0,\infty))$;
\item[\rm (vi)] $u(x,t)>0$ for every $x\in\mathbb{R}$, $t>0$.

\end{enumerate}
\end{Theorem}

{\sc Plan of the paper. } We will cover the existence and uniqueness
theory in sections \ref{sect-uniqueness} and \ref{sect-existence};
we borrow results and ideas from \cite{pqrv2}.  Section
\ref{sect-strong} is devoted to obtain some basic properties of the
solutions.

We then proceed with the smoothing effect, Section
\ref{sec.smoothing}, first from $L^1\cap L^p$, $p>1$, to $L^\infty$
and, then from $\cal X$ to $L^2$. The proof entails a number of new
ideas, in particular the use of a Trudinger inequality for
fractional exponents.

In Section \ref{sec.reg} we perform a delicate regularity
analysis to show that solutions are $C^\infty$ smooth in space and time,
and hence classical solutions of the equation.

We next describe in Section \ref{sec.trans} the transformation that
passes from the equation in~\eqref{eq:main} to the nonlocal
diffusion-transport model \eqref{eq:transport-rare}, and the results
obtained for the latter.

We finally include two appendixes.  The first one is
devoted to  a generalization of the Nash-Trudinger type inequality
used in the proof of the smoothing effect. In the second one we
consider another tool used in that proof, an interesting calculus
inequality.

\section{Uniqueness}\label{sect-uniqueness}
\setcounter{equation}{0}

As mentioned in the introduction, in order to prove uniqueness we
have to restrict the class of solutions under consideration. We will
give two results in this direction: in the first one we restrict
ourselves to weak solutions that satisfy $u\in
L^2(\mathbb{R}\times(0,T))$ for all $T>0$, and in the second to the
class of strong solutions.

\begin{Theorem}\label{th:uniqueness1}
Problem \eqref{eq:main} has at most one weak solution satisfying
$u\in L^2(\mathbb{R}\times(0,T))$ for all $T>0$.
\end{Theorem}
\noindent{\it Proof. } We adapt the classical uniqueness proof for
porous medium equations due to Oleinik, Kalashnikov and
Czou~\cite{OKC}.

Let $u$ and $\widetilde u$ be two weak
solutions to problem~\eqref{eq:main}. We subtract the weak
formulations for $u$ and $\widetilde u$ and take
$$
\varphi(x,t) =\left\{\begin{array}{ll}\displaystyle\int_t^T
(\log(1+u)-\log(1+\widetilde u))( x,s)\,ds,\qquad& 0\le t\le
T,\\[8pt]
0,\qquad& t\ge T,
\end{array}
\right.
$$
as a test function. Notice that, since the initial data of both
solutions coincide, we do not need $\varphi$ to vanish at $t=0$.
After an integration in time we get
$$
\begin{array}{l}
\displaystyle\int_0^T\int_{{\mathbb{R}}}(u-\widetilde
u)(x,t)(\log(1+u)-\log(1+\widetilde u))(x,t)\,dxdt
\\
[4mm] \qquad\qquad+\displaystyle\frac
12\int_{{\mathbb{R}}}\left(\int_0^T
(-\Delta)^{1/4}(\log(1+u)-\log(1+\widetilde u))(x,s)\,
ds\right)^2dx= 0.
\end{array}$$
The condition $u,\widetilde u\in L^2(\mathbb{R}\times(0,T))$ ensures
that the first integral is well defined. Since both integrands are
nonnegative, they must be identically zero. Therefore, $u=\widetilde
u$. \qed

\noindent\emph{Remark. } In particular,  for $f\in
L^1(\mathbb{R})\cap L^\infty(\mathbb{R})$ there is at most a bounded
weak solution.

To prove uniqueness in the class of strong solutions we use the
extension technique.

\begin{Theorem}\label{th:uniqueness2} If $u$ and $\widetilde u$ are strong solutions
to problem~\eqref{eq:main},  for every $0\le t_1<t_2$ we have
\begin{equation}
\label{eq:contractivity.uniqueness2}
\int_{{\mathbb{R}}}(u-\widetilde u)_+(x,t_2)\,dx\le
\int_{{\mathbb{R}}}(u-\widetilde u)_+(x,t_1)\,dx.
\end{equation}

\end{Theorem}

\noindent{\it Proof. } Let $p$ be a smooth monotone approximation to
the sign function such that $0\le p\le1$, and let $j$ be such that
$j'=p$, $j(0)=0$.  Let $\zeta\in C_0^\infty(\mathbb{R})$ be a
cut-off function, $0\le\zeta\le1$, $\zeta(x)=1$ for $|x|\le 1$,
$\zeta(x)=0$ for $|x|\ge2$, and $\zeta_R=\zeta(x/R)$.

Let $z=\log(1+u)-\log(1+\widetilde u)$. Using~\eqref{eq:plancherel}
and~\eqref{eq:fundamental.property.extensions} we get, for any
$0<t_1<t_2$,
$$
\begin{array}{l}
\displaystyle\int_{t_1}^{t_2}\int_{\mathbb{R}}\frac{\partial(u-\widetilde
u)}{\partial t}p(z)\zeta_R=
-\displaystyle\int_{t_1}^{t_2}\int_{\mathbb{R}}(-\Delta)^{1/4}z\,(-\Delta)^{1/4}\left(p(z)\zeta_R\right)\\[4mm]
\qquad\qquad=\displaystyle-\int_{t_1}^{t_2}\int_0^\infty\int_{\mathbb{R}}\nabla\E(z)\cdot\nabla\left(p(E(z))E(\zeta_R)\right)\\[4mm]
\qquad\qquad=\displaystyle
-\int_{t_1}^{t_2}\int_0^\infty\int_{\mathbb{R}}
(p'(\E(z))|\nabla\E(z)|^2E(\zeta_R)+ \nabla
j(\E(z))\cdot\nabla\E(\zeta_R))\\[4mm]
\qquad\qquad\le-\displaystyle\int_{t_1}^{t_2}\int_{\mathbb{R}}
(-\Delta)^{1/4}j(z)\,(-\Delta)^{1/4}\zeta_R\le\int_{t_1}^{t_2}\int_{\mathbb{R}}
j(z)|(-\Delta)^{1/2}\zeta_R|\\[4mm]
\qquad\qquad\le\displaystyle\frac
c{R}\int_{t_1}^{t_2}\int_{\mathbb{R}} |z| \le\frac
{c(t_2-t_1)}R\max_{t\in[t_1,t_2]}\max\{\|u(\cdot,t)\|_1,\|\widetilde
u(\cdot,t)\|_1\}.
\end{array}
$$
where we have used that
$|(-\Delta)^{1/2}\zeta_R(x)|=|(-\Delta)^{1/2}\zeta (x/R)|/R\le c/R$,
$0\le j(z)\le |z|$, and the fact that $\log(1+u)\le u$ for all $u\ge
0$. We end by letting $R\to\infty$ and $p$ tend to the sign
function. The case $t_1=0$ is obtained passing to the limit, using
 the $L^1$--continuity of $u(\cdot,t)$ at $t=0$.~\qed

\section{Existence of weak solutions}\label{sect-existence}
\setcounter{equation}{0}

The aim of this section is to construct a weak solution for any
initial data in $\mathcal{X}$. We will prove later, in
section~\ref{sect-strong}, that this solution, being strong, falls
within the uniqueness class.

\begin{Theorem}\label{th:plocal} For every $f\in\mathcal{X}$ there exists
a  weak solution $u$ to problem~\eqref{eq:main}. This solution
satisfies $u\ge0$,
\begin{equation}\label{energy}
\displaystyle\int_0^\infty\int_{\mathbb{R}}
|(-\Delta)^{1/4}\log(1+u)|^2\,dxdt  \le L_{\mathcal{X}}(f),
\end{equation}
and, if $f\in L^\infty(\mathbb{R})$,  $\|u(\cdot,t)\|_\infty\le
\|f\|_\infty$.
\end{Theorem}

\noindent{\it Proof. } The construction of the solution uses several
approximations. We refer to \cite{pqrv2} for the details,  where a
similar calculation is made for the fractional porous medium
equation~\eqref{eq:fpme}.

{\sc Step 1.} We first consider initial functions $f\in
L^1_+({\mathbb{R}})\cap L^\infty({\mathbb{R}})$. We  use the
formulation of the problem in the extension to $\mathbb{R}^2_+$
version \eqref{pp:local}. By means of the Crandall-Liggett Theorem
\cite{crandall-liggett} we are reduced to deal with the elliptic
related problem
\begin{equation}\label{eq:local-elliptic}
\left\{
\begin{array}{ll}
\Delta w=0,\qquad &x\in\mathbb{R},\,y>0,\\
-\partial_y w+\beta(w)=g,\qquad&x\in\mathbb{R},\,y=0,
\end{array}\right.
\end{equation}
with $g\in L^1_+({\mathbb{R}})\cap L^\infty({\mathbb{R}})$. Finally
we substitute the half space $\mathbb{R}^2_+$ by a half ball
$B_R^+=\{(x,y)\,:\,|x|^2+y^2<R^2,\,y>0\}$. We impose zero Dirichlet
data on the \lq\lq new part'' of the boundary. Therefore we are led
to study the problem
\begin{equation}
\left\{
\begin{array}{ll}
\Delta w=0\qquad &\mbox{in } B_R^+,\\
w=0\qquad &\mbox{on }
\partial B_R^+\cap\{y>0\},\\
-\partial_y w+\beta(w)=g\qquad&\mbox{on } D_R:=\{|x|<R,\, y=0\},
\end{array}
\right. \label{pp:local-elliptic-bdd}
\end{equation}
with $g\in L^\infty(D_R)$ given. Minimizing the functional
$$
J(w)=\frac12\int_{B_R^+}|\nabla
w|^2+\int_{D_R}\left(e^w-(1+g)w\right)
$$
in the admissible set $\mathcal{A}=\{w\in H^1(B_R^+): 0\le \beta(w)\le \|g\|_\infty\}$,
we obtain a unique solution $w=w_R$ to problem
\eqref{pp:local-elliptic-bdd}.
Moreover, if $g_1$ and $g_2$ are two admissible data, then the
corresponding weak solutions satisfy the $L^1$-contraction property
\begin{equation*}\label{eq:-contraction-bdd}
\int_{D_R}\left(\beta({w_1(x,0)})-\beta({w_2(x,0)})\right)_+\,dx\le
\int_{\mathbb{R}}\left(g_1(x)- g_2(x)\right)_+\,dx.
\end{equation*}

{\sc Step 2.} The passage to the limit $R\to\infty$ uses the
monotonicity in $R$ of the approximate solutions $w_R$. We obtain a
function $w_\infty=\lim_{R\to\infty}w_R$ which is a weak solution to
problem \eqref{eq:local-elliptic}. The above contractivity property
also holds in the limit. Moreover,
$\|\beta(w_\infty(\cdot,0))\|_{L^\infty(\mathbb{R})}\le\|g\|_{L^\infty(\mathbb{R})}$,
and  $w_\infty \ge0$, since $g\ge0$.

{\sc Step 3.} By the previous step, and using the Crandall-Liggett
Theorem, we obtain the existence of a unique mild solution
$\overline w$ to the evolution problem \eqref{pp:local}. To prove
that $\overline w$ is moreover a weak solution to problem
\eqref{pp:local}, one needs to show that it lies in the right energy
space. This is done using the same technique as in~\cite{pqrv},
which yields the energy estimate
\begin{equation*}\label{L2grad}
\int_0^T\int_0^\infty\int_{\mathbb{R}} |\nabla \overline
w(x,y,t)|^2\,dxdydt\le L_{\mathcal{X}}(f)\quad \text{for every }T>0.
\end{equation*}
Hence the function $u=\beta({\overline
w(\cdot,0)})$ is a  weak solution to problem \eqref{eq:main}. In addition, $\|\beta(\overline
w(\cdot,0))\|_{L^\infty(\mathbb{R}\times(0,\infty))}\le\|f\|_{L^\infty(\mathbb{R})}$,
and $\overline w \ge0$.  In
order to obtain estimate \eqref{energy} we recall the isometry
between  $\dot H^{1/2}(\mathbb{R})$ and $\mathcal{H}$. The Semigroup
Theory also guarantees that the constructed solutions satisfy the
$L^1$-contraction property $\|u(\cdot,t)-\widetilde
u(\cdot,t)\|_1\le \|f-\widetilde f\|_1$.

{\sc Step 4.} In this last step we consider general data
$f\in\mathcal{X}$. Let $\{f_k\}\subset L^1_+({\mathbb{R}})\cap
L^\infty({\mathbb{R}})$ be a sequence of functions converging to $f$
in $L^1(\mathbb{R})$,  and let $\{u_k\}$ be the sequence of the
corresponding solutions. Thanks to the $L^1$-contraction property we
know that $u_k(\cdot,t)\to u(\cdot,t)$ in $L^1({\mathbb{R}})$ for
all $t>0$ for some function $u$. Moreover, nonlinear Semigroup
Theory guarantees that $u_k\to u$ in $C([0,\infty):L^1(\mathbb{R}))$
\cite{Crandall}. On the other hand, using estimate~\eqref{energy},
we have $\log(1+u_k)\in L^2((\tau,\infty):\dot
H^{1/2}({\mathbb{R}}))$ uniformly in $k$. Thus the limit $u$ is a
weak solution to problem~\eqref{eq:main} for every $t\ge\tau$. The
$L^1$-contraction together with the $L^1$-continuity allow to go
down to $\tau=0$.~\qed

\section{Strong solutions and energy estimates}\label{sect-strong}
\setcounter{equation}{0}

We still have to prove that the weak solutions that we have
constructed are in fact strong. As a first step we consider the case
of bounded weak solutions. The general case will follow by
approximation as a consequence of the smoothing effect; see
Section~\ref{sec.smoothing}.

\begin{Proposition}\label{prop:estimates}
  Let $u$ be a bounded weak solution to problem  \eqref{eq:main}. Then $u$ is a strong solution and
  \begin{equation}\label{eq:u_t-en-L2}
    \int_t^\infty\int_{\mathbb{R}}|\partial_t u|^2\,dxds\le
    ct^{-1}(1+\|u(\cdot,t)\|_\infty)L_{\mathcal{X}}(f),\qquad t>0.
  \end{equation}
\end{Proposition}
\noindent{\it Proof. } In order to overcome the possible lack of
regularity in time,  we will work with the Steklov averages of
functions $g\in L^1_{\rm loc}(\mathbb{R}\times(0,\infty))$, defined
as
$$
g^h(x,t)=\frac1h\int_t^{t+h}g(x,s)\,ds.
$$
A similar approach is used for instance by B\'enilan and Gariepy in
\cite{Benilan-Gariepy} when dealing with evolution problems with
standard Laplacians. The use of Steklov averages makes  the process
rather technical. The estimates are simpler to obtain when we assume
regularity and work formally, and we invite the reader to do so.
However, such regularity cannot be assumed at this stage of the
theory.

Almost everywhere we have
$$
\partial_tg^h(x,t)=\delta^hg(x,t):=\frac{g(x,t+h)-g(x,t)}{h}.
$$
Let $h>0$. Given any $\varphi\in C^\infty_0(\mathbb{R}\times
(0,\infty))$, we may take $-\delta^{-h}\varphi$ as a test function
in the weak formulation. Then, using the \lq\lq integration by
parts'' formula $\int_0^\infty\int_\mathbb{R}\varphi\,\delta^h u =-
\int_0^\infty\int_\mathbb{R}u\delta^{-h}\varphi$,  we get that
$$
\int_0^\infty\int_{\mathbb{R}}\varphi\delta^h u\,dxdt=-\int_0^\infty\int_{\mathbb{R}}(-\Delta)^{1/4}
(\log(1+u))^h(-\Delta)^{1/4}\varphi\,dxdt.
$$
Taking  $\varphi=\zeta\partial_t(\log(1+u))^h$, where
$\zeta=\zeta(t)\in C_0^\infty((0,\infty))$, this identity becomes
\begin{equation}\label{eq:stek}
\begin{array}{rl}
\displaystyle\int_0^\infty\int_{\mathbb{R}}\zeta\partial_t
u^h\,\partial_t
(\log(1+u))^h\,dxdt&\displaystyle=-\frac12\int_0^\infty\int_{\mathbb{R}}\zeta\,
\partial_t\left|(-\Delta)^{1/4} (\log(1+u))^h\right|^2\,dxdt
\\ [3mm]&\displaystyle=\frac12\int_0^\infty\int_{\mathbb{R}}\zeta'\,
\left|(-\Delta)^{1/4} (\log(1+u))^h\right|^2\,dxdt.
\end{array}
\end{equation}
We now restrict ourselves to functions $\zeta$ which are cut-off
functions for the set $[t_1,t_2]$. To be more precise, we consider
$\psi\in C^\infty(\mathbb{R})$ such that $\psi'\ge0$, $\psi(t)=0$
for $t\le 1/2$, $\psi(t)=1$ for $t\ge1$, and then define
$\zeta(t)=\psi(t/t_1)-\psi(t/(2t_2))$.  Then, using that
$\zeta'(t)\le t_1^{-1}\max\psi'$, together with  the inequality
$\delta^h u\,\delta^h \log(1+u)\ge c\, (\delta^hu)^2$, (with
$c=(1+\|u\|_\infty)^{-1}$), we get
$$
c\int_{t_1}^{t_2}\int_{\mathbb{R}}(\delta^hu)^2\,dxdt\le
\frac1{2t_1}\int_0^\infty\int_{\mathbb{R}}\left|(-\Delta)^{1/4}
(\log(1+u))^h\right|^2\,dxdt.
$$
The energy estimate \eqref{energy} implies that the right-hand side
is bounded for $h$ small. Therefore there is a sequence $h_n\to0^+$
and a function $g\in L^2(\mathbb{R}\times(t,\infty))$ for all $t>0$
such that $\delta^{h_n}u\to g$ weakly in
$L^2(\mathbb{R}\times(t,\infty))$. It satisfies
$$
\int_{t_1}^{t_2}\int_{\mathbb{R}}g^2\,dxdt\le
ct_1^{-1}(1+\|u(\cdot,t_1)\|_\infty)\int_0^\infty\int_{\mathbb{R}}\left|(-\Delta)^{1/4}
\log(1+u)\right|^2\,dxdt.
$$

On the other hand,
$$
\begin{aligned}
-\int_0^\infty\int_\mathbb{R}u\partial_t\varphi\,dxdt&=-\lim_{h_n\to0^+}\int_0^\infty\int_\mathbb{R}u\delta^{-h_n}\varphi\,dxdt\\
&=
\lim_{h_n\to0^+}\int_0^\infty\int_\mathbb{R}\delta^{h_n}u\varphi\,dxdt\\
&=\int_0^\infty\int_\mathbb{R}g\varphi\,dxdt,
\end{aligned}
$$
which means that the distributional derivative $\partial_tu$ is in
fact a function that coincides with $g$ almost everywhere. \qed

We next prove that the $L^p$-norms do not increase with time. The
main tool, used also later in the proof of the smoothing effect,
Section 6, is the generalized Stroock-Varopoulos inequality
\cite{Stroock-1984}, \cite{Varopoulos-1985},
\begin{equation}
\label{eq:strook.varopoulos2} \int_{\mathbb{R}}A(z)(-\Delta)^{1/2} z
\ge \int_{\mathbb{R}}\left|(-\Delta)^{1/4}B(z)\right|^2,
\end{equation}
where $A'=(B')^2$. An easy proof using the local realization of the
half-Laplacian  (in a more general setting)  is given in~\cite[Lemma
5.2]{pqrv2}.

\begin{Proposition}\label{prop:decay-Lp}
Let $u$ be a bounded weak solution to problem  \eqref{eq:main}.
Then, for every  $0\le t_1<t_2$ we have
$$
L_{\mathcal{X}}(u(\cdot,t_2))\le
L_{\mathcal{X}}(u(\cdot,t_1)),\qquad
\|u(\cdot,t_2)\|_p\le\|u(\cdot,t_1)\|_p, \quad 1\le p\le\infty.
$$
\end{Proposition}
\noindent{\it Proof. }  The first estimate is obtained directly
multiplying the equation by $\log(1+u)$, as mentioned in
Section~\ref{sect-main.results}. The cases $p=1$ and $p=\infty$ in
the second inequality follow from the elliptic estimates in
Section~\ref{sect-existence}. For the rest of the cases, we  put
$A(z)=u^{p-1}$, $z=\log(1+u)$ in \eqref{eq:strook.varopoulos2}.
Since $(-\Delta)^{1/2} z\in L^2(\mathbb{R})$ a.e.~in $t$, if
$p\ge3/2$ we have $A(z)\in L^2(\mathbb{R})$. Assume this is the
case. We then multiply the equation by $A(z)$ and integrate in
$\mathbb{R}\times(t_1,t_2)$ to obtain
$$
\frac1p\int_{\mathbb{R}}\Big(u^p(x,t_2)-u^p(x,t_1)\Big)\,dx\le
-\int_{t_1}^{t_2}\int_{\mathbb{R}}
\left|(-\Delta)^{1/4}G(u)(x,t)\right|^2\,dxdt\le0,
$$
where $G(u)=B(z)=\int_0^{u}\sqrt{(p-1)s^{p-2}/(1+s)}\,ds$.

For the case $1<p<3/2$, we approximate the function $A(z)$ by
$$
A_\varepsilon(z)=\left\{\begin{array}{ll} u^{p-1}&\quad\mbox{for }
u\ge\varepsilon, \\
\varepsilon^{p-2}u&\quad\mbox{for } 0\le u<\varepsilon,
\end{array}\right.$$
and then let $\varepsilon$ tend to zero. \qed

The $L^1$-norm is not only non-increasing; it is
conserved.\begin{Theorem}\label{th:mass} Let $u$ be a strong
solution to problem~\eqref{eq:main}. For every $t>0$ we have
$$
\int_{{\mathbb{R}}}u(x,t)\,dx= \int_{{\mathbb{R}}}f(x)\,dx.
$$
\end{Theorem}
\noindent{\it Proof.}  We take a nonnegative non-increasing cut-off
function $\psi(s)$ such that $\psi(s)=1$ for $0\le s\le1$,
$\psi(s)=0$ for $s\ge2$, and define $\phi_R(x)=\psi(|x|/R)$. Observe
that
$|(-\Delta)^{1/2}\phi_R(x)|=R^{-1}|(-\Delta)^{1/2}\psi(|x|/R)|\le
c/R$. Multiplying the equation by $\phi_R$  and integrating by
parts, we obtain, for every $t_2>t_1>0$,
$$\begin{array}{rl}
\displaystyle \left|\int_{\mathbb{R}}\Big(
u(x,t_2)-u(x,t_1)\Big)\phi_R(x)\,dx\right|
&\displaystyle=\left|\int_{t_1}^{t_2}\int_{{\mathbb{R}}}\log(1+u)(x,t)\,
(-\Delta)^{1/2}\phi_R(x)\,dxdt\right| \\ [4mm]&\displaystyle\le c
R^{-1}  \max_{t\in[t_1,t_2]}\|u(\cdot,t)\|_1.
\end{array}
$$
In the last step we have used that $\log(1+u)\le u$. The result is
then obtained just passing to the limit $R\to\infty$. \qed

 Weak bounded solutions turn out to have an energy
which is well defined for all positive times.
\begin{Proposition}\label{prop-energy-dec}
Let $u$ be a bounded weak solution to problem \eqref{eq:main}. The
energy
\begin{equation*}
E(t)=\frac12\int_{\mathbb{R}}|(-\Delta)^{1/4}\log(1+u)(x,t)|^2\,dx
\end{equation*}
is a continuous function  in $(0,\infty)$ which does not increase
with time. Moreover,
\begin{equation}\label{eq:energy.estimate}
    E(t)\le
    (2t)^{-1} L_{\mathcal{X}}(f)\quad\text{for every }t>0.
\end{equation}

\end{Proposition}
\noindent{\it Proof. } Passing to the limit $h\to0$ in the
 identity \eqref{eq:stek}, we
 get
$$
\int_0^\infty\int_{\mathbb{R}}\zeta\frac{|\partial_t
u|^2}{1+u}\,dxdt=\frac12\int_0^\infty\int_{\mathbb{R}}\zeta'\,
\left|(-\Delta)^{1/4} \log(1+u)\right|^2\,dxdt
$$
for any test function $\zeta\in C_0^\infty((0,\infty))$. This means
that, as a distribution, $E'$ coincides with the function
$-\int_{\mathbb{R}}\frac{|\partial_t u|^2}{1+u}\,dx\le0$. Since the
latter belongs to $L^1(\mathbb{R})$, we conclude that $E\in
W^{1,1}((0,\infty))$, and therefore that it is a continuous
function. Now we have $ L_{\mathcal{X}}(f)\ge2 \int_{0}^{t}
E(s)\,ds\ge2 tE(t)$. \qed

In addition to the homogeneous Sobolev space $\dot H^{1/2}(\mathbb{R})$,
the function $\log(1+u)(\cdot,t)$ also belongs to the full space $H^{1/2}(\mathbb{R})$.

\begin{Proposition}\label{prop:H12}
Let $u$ be a bounded weak solution to problem \eqref{eq:main}.
Then  for every $t>0$
\begin{equation}\label{eq:Phi(u)-en-H1/2}
\|\log(1+u)(\cdot,t)\|_{H^{1/2}} \le
t^{-1/2}\big(L_{\mathcal{X}}(f)\big)^{1/2}+
 ct^{-1/4}\|f\|_1^{1/2}\big(L_{\mathcal{X}}(f)\big)^{1/4}.
\end{equation}
\end{Proposition}
\noindent{\it Proof. } Let $w=\log(1+u)$. We use interpolation and
the Nash-Gagliardo-Nirenberg inequality
\eqref{eq:gagliardo.nirenberg.type.inequality2} with $N=1$,
$\gamma=1/2$, $q=2$, $p=1$,  to get
$$
\|w(\cdot,t)\|_2\le\|w(\cdot,t)\|_3^{3/4}\|w(\cdot,t)\|_1^{1/4} \le
c\|(-\Delta)^{1/4}w(\cdot,t)\|_2^{1/2}\|w(\cdot,t)\|_1^{1/2}.
$$
Next we use that $\log(1+u)\le u$, the energy estimate
\eqref{eq:energy.estimate} and the conservation of mass to conclude
that
\begin{equation}
\label{eq:L2.estimate.log}
\|\log(1+u)(\cdot,t)\|_2\le ct^{-1/4}\big(L_{\mathcal{X}}(f)\big)^{1/4}\|u(\cdot,t)\|_1^{1/2}=
ct^{-1/4}\big(L_{\mathcal{X}}(f)\big)^{1/4}\|f\|_1^{1/2}.
\end{equation}
\qed

To end this section, we improve the regularity of $u$ and
$\log(1+u)$, giving an $L^2$-control of their gradients.

\begin{Corollary}Let $u$ be a bounded weak solution to problem \eqref{eq:main}. Then  $u$ and $\log(1+u)$ belong to $
L^2_{\rm loc}((0,\infty):H^1(\mathbb{R}))$, and
$$
\begin{array}{l}
\displaystyle
\int_t^\infty\int_{\mathbb{R}}|\partial_x\log(1+u)(x,s)|^2\,dxds\le
ct^{-1}(1+\|u(\cdot,t)\|_\infty)L_{\mathcal{X}}(f), \\ [3mm]
\displaystyle
\int_t^\infty\int_{\mathbb{R}}|\partial_xu(x,s)|^2\,dxds\le
ct^{-1}(1+\|u(\cdot,t)\|_\infty)^3L_{\mathcal{X}}(f).
\end{array}$$
\end{Corollary}
\noindent{\it Proof. } It is clear from~\eqref{eq:L2.estimate.log}
that  $\log(1+u)\in L^2((0,T)\,:\,L^2(\mathbb{R}))$ for every $T>0$.
To estimate the gradient we just use \eqref{eq:u_t-en-L2} and the
equation. Actually,
$$
\begin{array}{l}
\displaystyle\int_t^\infty\int_{\mathbb{R}}|\partial_x\log(1+u)(x,s)|^2\,dxds
=\int_t^\infty\int_{\mathbb{R}}|(-\Delta)^{1/2}\log(1+u)(x,s)|^2\,dxds
\\ [4mm]
\qquad=\displaystyle\int_t^\infty\int_{\mathbb{R}}|\partial_tu(x,s)|^2\,dxds
\le   c t^{-1}(1+\|u(\cdot,t)\|_\infty)L_{\mathcal{X}}(f).
\end{array}
$$
As to $u$, we just observe that $\partial_x
u=(1+u)\partial_x\log(1+u)$. \qed

\section{Smoothing effect} \label{sec.smoothing}
\setcounter{equation}{0}

In Section \ref{sect-existence} we have constructed a weak solution
of problem~\eqref{eq:main} for general initial data
$f\in\mathcal{X}$ by approximation with initial data in
$L^1(\mathbb{R})\cap L^\infty(\mathbb{R})$. Our next aim is to prove
that this solution becomes immediately bounded; in particular it is
strong. Boundedness will follow from an estimate for bounded weak
solutions, formula~\eqref{eq:full.smoothing.effect}, which
\emph{does not depend} on the $L^\infty$ norm of the datum, but only
on $L_{\mathcal{X}}(f)$ (and time).

The result will be obtained by combining  $L^2\to L^\infty$ and
$\mathcal{X}\to L^2$ smoothing effects. The  $L^2\to L^\infty$
result is in fact a particular instance of a more general $L^p\to
L^\infty$ result, valid for all $p>1$.

\begin{Theorem}\label{th:smoothing1} Let $u$ be a bounded weak solution,
and let $p>1$. There is a constant $C>0$ that depends only on $p$
such that
\begin{equation} \|u(\cdot,t)\|_\infty\le C\,
\max\{t^{-1/(p-1) }\|f\|_{p}^{p/(p-1)},t^{-1/p }\|f\|_{p}\}.
\label{eq:L-inf}\end{equation}
\end{Theorem}

We recall that the corresponding formula for the fractional PME with
$m>0$ reads, in the case $N=\sigma=1$,
\begin{equation}
\|u(\cdot,t)\|_\infty\le C\,t^{-1/(m+p-1) }\|f\|_{p}^{p/(m+p-1)},
\label{smoothing-m>0}\end{equation} for every $p\ge1$, cf.
\cite{pqrv2}. Observe that when $m=0$ these exponents make  sense
for $p>1$ but not for $p=1$.  It is also worth noticing that formula
\eqref{eq:L-inf} can be obtained by formally putting in
\eqref{smoothing-m>0} $m=0$ for $u$ large and $m=1$ for $u$ small.

\noindent{\it Proof.} The proof follows the same Moser iterative
technique used in \cite{pqrv2}, but it is a little more involved.
Let $t>0$ be fixed, and consider the sequence of times
$t_k=(1-2^{-k})t$, $p_k=2^kp$. We multiply the equation (recall that
it is satisfied a.e. since $u$ is a strong solution) by the test
function
$$
\phi=\frac{u^{p_k-1}}{p_k-1}+\frac{u^{p_k}}{p_k}
$$
and integrate in $\mathbb{R}\times(t_k,t_{k+1})$ (for
$p_0=p\in(1,3/2)$ we need an extra approximation argument,  as in
Proposition~\ref{prop:decay-Lp}, to justify the computation). Using
now the Stroock-Varopoulos inequality \eqref{eq:strook.varopoulos2},
we get
$$
\begin{array}{rl}
\displaystyle\frac1{p_k(p_k-1)}\|u(\cdot,t_k)\|_{p_k}^{p_k}
&\displaystyle+\frac1{p_k(p_k+1)}\|u(\cdot,t_k)\|_{p_k+1}^{p_k+1} \\
[3mm] &\displaystyle\ge
\frac4{p_k^2}\int_{t_k}^{t_{k+1}}\int_{\mathbb{R}}|(-\Delta)^{1/4}u^{p_k/2}(x,\tau)|^2\,dxd\tau.
\end{array}
$$
Multiplying and dividing by  $\|u(\cdot,\tau)\|_{r}^{r}$, for some
$r>1$, $r\ge p_k/2$, using that the $L^r$ norms do not increase in
time, and applying  the Nash-Gagliardo-Nirenberg type inequality
\eqref{eq:gagliardo.nirenberg.type.inequality2} with $N=1$,
$\gamma=1/2$, $q=2$, we get
\begin{equation}
\label{eq:recursion}
\|u(\cdot,t_{k+1})\|_{p_k+r}^{p_k+r}\le c2^kt^{-1}
\|u(\cdot,t_k)\|_{r}^{r}\Big(\|u(\cdot,t_k)\|_{p_k}^{p_k}+\|u(\cdot,t_k)\|_{p_k+1}^{p_k+1}\Big).
\end{equation}
Let us denote
$U_k=\max\{\|u(\cdot,t_k)\|_{p_k},\,\|u(\cdot,t_k)\|_{p_k+1}^{(p_k+1)/p_k}\}$.
Taking $r=p_k$ and $r=p_k+1$ in \eqref{eq:recursion} we get that
both $\|u(\cdot,t_{k+1})\|_{p_{k+1}}^{p_{k+1}}$ and
$\|u(\cdot,t_{k+1})\|_{p_{k+1}+1}^{p_{k+1}+1}$ are smaller than
$c2^kt^{-1}U_k^{p_{k+1}}$, from where we obtain
$$
U_{k+1}\le (c2^kt^{-1})^{1/(p_{k+1})}U_k= (c2^{
k/p}t^{-1/p})^{1/2^{k+1}}U_k.
$$
This recursive relation yields
$$
\|u(\cdot,t)\|_\infty=\lim_{k\to\infty}U_k\le ct^{-1/p}U_0
=ct^{-1/p}\max\{\|f\|_{p}\,,\|f\|_{p+1}^{(p+1)/p}\}.
$$

The final step is to get rid of the $L^{p+1}$-norm. Using H\"older's
inequality and the decay of the $L^p$-norms we get
$$
\|u(\cdot,t)\|_\infty\le
c(t/2)^{-1/p}\|f\|_{p}\max\{1,\|u(\cdot,t/2)\|_\infty^{1/p}\}.
$$
If $\|u(\cdot,t/2)\|_\infty\le1$, then $\|u(\cdot,t)\|_\infty\le
c2^{1/p}t^{-1/p}\|f\|_{p}$ and we are done.
If, on the contrary,  $\|u(\cdot,t/2)\|_\infty\ge1$,  we have
$$
\|u(\cdot,t)\|_\infty\le
c(t/2)^{-1/p}\|f\|_{p}\|u(\cdot,t/2)\|_\infty^{1/p}.
$$
Since in this case, by the maximum principle, we have
$\|u(\cdot,\tau)\|_\infty\ge1$ for every $0<\tau<t/2$, we may
iterate this estimate to get
\begin{equation*}\label{eq:smooth1}
\|u(\cdot,t)\|_\infty\le ct^{-1/(p-1)}\|f\|_{p}^{p/(p-1)}.
\end{equation*}
 \qed

The above method does not allow to go down to $p=1$. This drawback
was already present in the PME case (both local and nonlocal, see
\cite{vazquez} and \cite{pqrv2}), where the limit exponent was
$p=\max\{1,\,(1-m)N/\sigma\}$. In the case $N=\sigma=1$, and putting
$m=0$, we get that the limit exponent should be $p=1$, but it is not
clear if solutions will become bounded when  the initial datum only
belongs to $L^1(\mathbb{R})$. Nevertheless, we may consider initial
values in the slightly smaller space ${\cal X}$. This is our next
goal.

\begin{Theorem}\label{th:smoothing2} Let $u$ be a bounded weak solution. There is a constant $C>0$ such that
\begin{equation}\label{eq:u-en-L2}
\int_{\mathbb{R}}u^2(x,t)\,dx\le \mbox{\rm exp}\left\{C\left(t^{-1/2}\big(L_{\mathcal{X}}(f)\big)^{1/2}+
t^{-1/4}\|f\|_1^{1/2}\big(L_{\mathcal{X}}(f)\big)^{1/4}\right)\right\}-1.
\end{equation}
\end{Theorem}

\noindent{\it Proof. } Fix any time $t>0$ and let
$w=\log(1+u(\cdot,t))$. We know from Proposition~\ref{prop:H12} that $w\in H^{1/2}(\mathbb{R})$. Hence, using
the Trudinger type inequality \eqref{eq:Strichartz},   with $N=1$
and $\gamma=1/2$, we obtain
$$
\int_{\mathbb{R}}\Big(e^{w^2/c\|w\|^2_{H^{1/2}}}-1\Big)\le 1.
$$
We now apply the  calculus inequality $ (e^{w}-1)^2\le
(e^k-1)(e^{w^2/k}-1) $ (see Lemma~\ref{lem:exp-k} below  for the
proof), to get
\begin{equation*}
\int_{\mathbb{R}}u^2(x,t)\,dx\le
\mbox{exp}(c\|w\|^2_{H^{1/2}})-1.
\end{equation*}
We conclude using the energy estimate \eqref{eq:Phi(u)-en-H1/2}.\qed

To obtain the $\mathcal{X}-L^\infty$ smoothing effect we just have to combine  Theorems~\ref{th:smoothing1} and~\ref{th:smoothing2}.

\

\noindent\emph{Proof of Theorem~\ref{th:props}\text{\rm -(ii)}.}  We
first consider the case of initial data which are moreover bounded.
The general case is dealt with by approximation.

Using the $L^p-L^\infty$ estimate  \eqref{eq:L-inf} with $p=2$, and
the $\mathcal{X}$-$L^2$  estimate \eqref{eq:u-en-L2} we get, first
for $t$ small,
$$
\|u(\cdot,t)\|_\infty\le Ct^{-1}\|u(\cdot,t/2)\|_2^2\le
Ct^{-1}\mbox{exp}(Ct^{-1/2}L_{\mathcal{X}}(f)^{1/2}),
$$
and then for $t$ large
$$
\|u(\cdot,t)\|_\infty\le Ct^{-1/2}\|u(\cdot,t/2)\|_2\le C
t^{-3/4}\|f\|_1^{1/2}L_{\mathcal{X}}(f)^{1/4}.
$$
Combining both estimates we obtain \eqref{eq:full.smoothing.effect}.
\qed

\section{Regularity and positivity}\label{sec.reg}
\setcounter{equation}{0}

The  solution that we have constructed in the previous sections is
$C^\infty$ for all positive times, and hence classical. This is the
content of the present section.

\subsection{$C^{1,\alpha}$ regularity} The first and more difficult step is to prove that the solution $u$ is
$C^{1,\alpha}$ for all $\alpha\in(0,1)$. Actually, given $\tau>0$, $u$ is uniformly $ C^{1,\alpha}$ in $Q_\tau=\mathbb{R}\times(\tau,\infty)$, denoted $u\in C^{1,\alpha}_{\rm u}(Q_\tau)$ for short.

\begin{Theorem}
\label{th:regularity} Let $f\in\mathcal{X}$.  The strong  solution
to  problem \eqref{eq:main} satisfies $u\in C^{1,\alpha}_{\rm u}(Q_\tau)$ for every $0<\alpha<1$
and $\tau>0$.
 \end{Theorem}

\noindent{\it Proof. }  {\sc Step 1:} \emph{ $u\in
C^\alpha_{\rm u}(Q_\tau)$ for some $\alpha\in(0,1)$ and every $\tau>0$.}

\noindent Once we know that $u$ is bounded in the time interval $t\ge\tau>0$,
the result follows from the regularity results for problem
\eqref{pp:local} from
Athanasopoulos-Caffarelli \cite{ac}, since the nonlinearity $\beta(u)$ satisfies the
non-degeneracy condition required in that paper.

\medskip

\noindent {\sc Step 2:} \emph{$u\in C^\alpha_{\rm u}(Q_\tau)$
for every $0<\alpha<1$ and every $\tau>0$.}

\noindent To prove this, we will show that H\"older regularity can
be \lq\lq doubled'', following ideas from Caffarelli and Vasseur
\cite{cv}; i.e., if $u\in C^{\alpha}_{\rm u}(Q_\tau)$ for some
$\alpha\in(0,1/2)$, then $u\in C^{2\alpha}_{\rm u}(Q_\tau)$. The
claimed regularity is then obtained repeating the argument a finite
number of times.

Let $(x_0,t_0)\in Q_\tau$   be fixed, and denote $u_0=u(x_0,t_0)$.
We write the equation in \eqref{eq:main} as a fractional linear heat
equation with a (nonlinear) source term,
\begin{equation}
\label{eq:semilinear}
\partial_tu+\mu(-\Delta)^{1/2} u=-(-\Delta)^{1/2} (\log(1+u)-\mu u).
\end{equation}
If we take $\mu=1/(1+u_0)$, the right-hand side of equation
\eqref{eq:semilinear}  can be written as $-(-\Delta)^{1/2}F(u)$,
where $F(u)=\log(\mu(1+u))-\mu(u-u_0)$ satisfies $F(u_0)=F'(u_0)=0$.
After a time shift, we may assume that $u$ is uniformly $C^\alpha$
and bounded down to $t=0$. Recall now that the fundamental solution
to the fractional heat equation $\partial_t u+(-\Delta)^{1/2}u=0$ is
the Poisson kernel
$$
P(x,t)=\dfrac1\pi\dfrac{t}{x^2+t^2}.
$$
Taking a smooth approximation of $P(x,\mu t)$ as a test function in
the distributional version of~\eqref{eq:semilinear}, and passing to
the limit in the approximation we get that the solution $u$ can be
represented in the (mild solution) form
\begin{equation}\label{var-param}
\begin{array}{rl}
u(x,t)&\displaystyle=\int_{\mathbb{R}}P(x-x_1,\mu t)f(x_1)\,dx_1
\\ [4mm]
&\displaystyle-\int_0^t\int_{\mathbb{R}}(-\Delta)^{1/2}P(x-x_1,\mu
(t-t_1))F(u(x_1,t_1))\,dx_1dt_1.
\end{array}
\end{equation}
The first term in the right-hand side of~\eqref{var-param} is regular, so we
concentrate on the second one.

We will use the notation  $y=(x,t)$ for the space-time variable, and
also $\overline y=(x,\mu t)$ to accommodate the distortion in time
created by the factor $\mu$. We are thus led to study the regularity
for the function
\begin{equation}
  \label{eq:g}
  g(y)=\int_{\mathbb{R}_+^2}A(\overline y-\overline y_1)\chi_{\{t_1<t\}}F(u(y_1))\,dy_1,
\end{equation}
where
$$
A(y)=A(x,t)\equiv(-\Delta)^{1/2}
P(x,t)=\dfrac1\pi\dfrac{x^2-t^2}{(x^2+t^2)^2}.
$$

Let us see first that the function $g$ is well defined. To this aim
we decompose $\mathbb{R}^2_+$ as $E_\rho\cup E^c_\rho$, where
$E_\rho$ is the ellipse
$$ E=\{y_1\in
\mathbb{R}_+^2\,:\,|\overline y_1-\overline y|<\rho\},
$$
with $\rho$  small. Observing that
$$
\int_{E_\rho}A(\overline y-\overline
y_1)\chi_{\{t_1<t\}}\,dy_1=\frac1{\pi\mu}\int_{\{x^2+t^2<\rho^2,\,t<0\}}\frac{x^2-t^2}{(x^2+t^2)^2}\,dxdt=0,
$$
we may write
$$
\begin{array}{rcl}
\displaystyle\left|\int_{E_\rho}A(\overline y-\overline
y_1)\chi_{\{t_1<t\}}F(u(y_1))\,dy_1\right|&\le&\displaystyle\int_{E_\rho}|A(\overline
y-\overline y_1)|\,|F(u(y_1))-F(u(y))|\,dy_1 \\ [4mm]
&\le&\displaystyle c\int_{E_\rho}\frac{dy_1}{|y-y_1|^{2-\beta}}\le c
\end{array}$$
for some $\beta>0$, since both $u$ and $F$ are H\"older continuous
functions. On the other hand, using that $F(u)$ is bounded we have
$$
\left|\int_{E^c_\rho}A(\overline y-\overline
y_1)\chi_{\{t_1<t\}}F(u(y_1))\,dy_1\right|\le
c\int_0^t\int_{|x_1-x|>1/(2\mu)}\frac{dx_1dt_1}{|x-x_1|^2}\le c.
$$

Now we will see that $g(y)$ has the same regularity as $F(u(y))$.
The key point is that if $u$ is $C^\alpha$ at $y_0$, then $F(u(y))$
is $C^{2\alpha}$ at $y_0$. Indeed, since $F(u_0)=F'(u_0)=0$, and
$|F''(u)|\le c$ (recall that $u\ge0$),  we have
$$
|F(u(y))|\le c|u(y)-u_0|^2\le c|y-y_0|^{2\alpha}
$$
for every $y\in \mathbb{R}_+^2$. Moreover, the constants are
independent of the point $y_0$. We observe also that $|A(\overline
y)|\le c|A(y)|$, where $c=c(\mu)$. Since $u$ is bounded and
nonnegative, the constant $c(\mu)$ can be taken independent of
$\mu$.

\begin{figure}[ht]
\begin{center}
\psfig{file=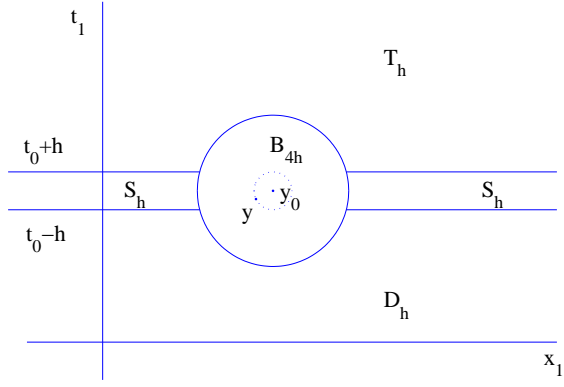,width=7.5cm}
\end{center}
\caption{Integration regions.}
\label{fig:integration.regions}
\end{figure}

Let $y\in \mathbb{R}_+^2$ be any point with $|y-
y_0|=h$. We have to prove that  the difference
\begin{equation}\label{eq:g-g0}
g(y_0)-g(y)=\int_{\mathbb{R}_+^2}\Big(A(\overline
y_0-\overline y_1)\chi_{\{t_1<t_0\}}-A(\overline y-\overline y_1)\chi_{\{t_1<t\}}\Big)F(u(y_1))\,dy_1
\end{equation}
is $O(h^{2\alpha})$ for $h$ small. In order to estimate the integral
in \eqref{eq:g-g0} we decompose $\mathbb{R}^2_+$ into four regions,
depending on the sizes of $|x_1-x_0|$ and $t_1-t_0$, see
Figure~\ref{fig:integration.regions}.

\begin{figure}[ht]
\begin{center}
\psfig{file=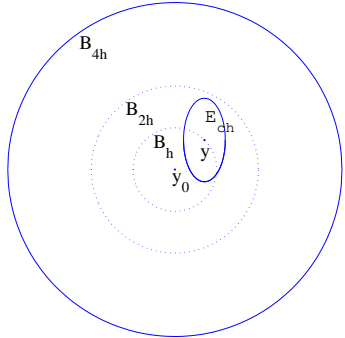,width=4.5cm}
\end{center}
\caption{Integration subregions in $B_{4h}$.}
\label{fig:integration.subregions}\end{figure}

\noindent (i) \emph{The small ball $
B_{4h}=\{|y_1-y_0|<4h\}\subset\mathbb{R}_+^2$. } The difficulty in
this region is the non-integrable singularity of $A(\overline y)$ at
$\overline y=0$. Integrability will be gained thanks to the
regularity of $F(u)$. We have,
$$
\int_{B_{4h}}|A(\overline y_0-\overline y_1)|\,|F(u(y_1))|\,dy_1\le
c\int_{B_{4h}}\frac{dy_1}{|y_1-y_0|^{2-2\alpha}}\le ch^{2\alpha}.
$$
In order to estimate $\int_{B_{4h}}A(\overline y-\overline
y_1)F(u(y_1))\,dy_1$, we consider as before the ellipse $ E_{ch}$,
where $c=c(\mu)$ is chosen to have $E_{ch}\subset B_{2h}$, see
Figure~\ref{fig:integration.subregions}.

We get
$$
\begin{array}{l}
\displaystyle\int_{B_{4h}}A(\overline y-\overline
y_1)\chi_{\{t_1<t\}}F(u(y_1))\,dy_1=
\\[4mm]
\qquad\qquad\underbrace{\displaystyle\int_{B_{4h}}A(\overline
y-\overline
y_1)\chi_{\{t_1<t\}}\Big(F(u(y_1))-F(u(y))\Big)\,dy_1}_{I_1}\\[4mm]
\qquad\qquad+ \underbrace{F(u(y))\int_{B_{4h}-E_{ch}}A(\overline
y-\overline y_1)\chi_{\{t_1<t\}}\,dy_1}_{I_2},
\end{array}
$$
since the integral over the ellipse is again zero by symmetry. To
estimate $I_1$ we use the Mean Value Theorem applied to the function
$F$ to see that
$$
|F(u(y))-F(u(y_1))|=|F'(\theta)|\,|u(y)-u(y_1)|\le
c\max\{|u(y)-u_0|,|u(y_1)-u_0|\}|y-y_1|^{\alpha},
$$
where $\theta$ is some value between $u(y)$ and $u(y_1)$. Therefore,
$$
|I_1|\le c\int_{B_{4h}}\frac1{|y-y_1|^2}|y-y_1|^\alpha
(|y_1-y_0|^\alpha+|y-y_0|^\alpha)\,dy_1\le ch^{2\alpha}.
$$
As to $I_2$, since we are far from the singularity of $A$,
$$
|I_2|\le c h^{2\alpha}\int_{B_{4h}-E_{ch}}\frac{dy_1}{h^2}\le
ch^{2\alpha}.
$$

\noindent (ii) \emph{The narrow strip
$S_h=\{|y_1-y_0|>4h,\;|t_1-t_0|<h\}$. } In this region we have
$|y_0-y_1|\le \frac43|y-y_1|$ and $|x_1-x_0|>3h$. Therefore,
$$
\begin{array}{l}
\displaystyle\int_{S_h}
|A(\overline y_0-\overline y_1)\chi_{\{t_1<t_0\}}-A(\overline y-\overline y_1)\chi_{\{t_1<t\}}||F(u(y_1))|\,dy_1\\[4mm]
\qquad\displaystyle\le \int_{S_h}\Big(|A(\overline y_0-\overline
y_1)|+|A(\overline y-\overline y_1)|\Big)|F(u(y_1))|\,dy_1 \le
\int_{S_h}\frac{dy_1}{|y_0-y_1|^{2-2\alpha}}\\[4mm]
\qquad\displaystyle\le
c\int_{t_0-h}^{t_0+h}\int_{|x_1-x_0|>3h}\frac{dx_1dt_1}{|x_0-x_1|^{2-2\alpha}}
\displaystyle\le ch^{2\alpha}.
\end{array}
$$

\noindent (iii) \emph{The complement of the ball $B_{4h}$ for large times, $T_h=\{|y_1-y_0|>4h,\;t_1>t_0+h\}$. } The integral
in this region is 0, since here we have
$$
A(\overline y_0-\overline y_1)\chi_{\{t_1<t_0\}}=A(\overline
y-\overline y_1)\chi_{\{t_1<t\}}=0.
$$

\medskip

\noindent (iv) \emph{The complement of the ball $B_{4h}$ for small
times, $D_h=\{|y_1-y_0|>4h,\;t_1<t_0-h\}$. } The required estimate
is obtained here using the fact that we are  integrating  a
difference of $A$'s, so there will be some cancelation. Indeed, by
the Mean Value Theorem,
$$
|A(\overline y_0-\overline y_1)-A(\overline y-\overline y_1)|\le
|\overline y_0-\overline y| \max\{|\partial_x A(\xi)|, |\partial_t
A(\xi)|\}\le ch/|\xi|^3,
$$
where $\xi=s(\overline y_0-\overline y_1)+(1-s)(\overline
y-\overline y_1)$ for some $s\in(0,1)$. On the other hand, since we
are in $D_h$, $|\overline y_0-\overline
y_1|\le\mu^{1/2}|y_0-y_1|\le\frac{4\mu(1-s)}{3} |\xi|\le c|\xi|$,
and we conclude that
$$
|A(\overline y_0-\overline y_1)\chi_{\{t_1<t_0\}}-A(\overline
y-\overline y_1)\chi_{\{t_1<t\}}|\le  \frac{ch}{|y_0-y_1|^3}.
$$
Therefore, assuming that $\alpha<1/2$,
$$
\begin{array}{l}
\displaystyle\int_{D_h}|A(\overline y_0-\overline y_1)\chi_{\{t_1<t_0\}}-A(\overline y-\overline y_1)\chi_{\{t_1<t\}}|\,|F(u(y_1))|\,dy_1 \\
\hspace{25mm}\displaystyle\le
ch\int_{D_h}\frac{dy_1}{|y_0-y_1|^{3-2\alpha}}\le ch^{2\alpha}.
\end{array}
$$

\medskip

\noindent {\sc Step 3:} \emph{$u\in
C_{\rm u}^{1,\alpha}(Q_\tau)$ for every $0<\alpha<1$ and every $\tau>0$.}

\noindent  We may assume, after a time shift, that $\tau=0$. Let
$z=y-y_0$. The result will follow from an estimate of the quantity
$$
\displaystyle
g(y_0+z)-2g(y_0)+g(y_0-z)=\int_{\mathbb{R}_+^2}\mathcal{A}(y_0,y,y_1)F(u(y_1))\,dy_1,
$$
where
$$
\displaystyle\mathcal{A}(y_0,y,y_1)=A(\overline y-\overline
y_1)\chi_{\{t_1<t\}}-2A(\overline y_0-\overline
y_1)\chi_{\{t_1<t_0\}}+A(2\overline y_0-\overline y-\overline
y_1)\chi_{\{t_1<2t_0-t\}}.
$$
As in the previous step, we consider separately the contributions to
the integral of the four regions shown in
Figure~\ref{fig:integration.regions}. The contribution of the ball
$B_{4h}$ is decomposed as the sum $J_1-2J_2+J_3$, where
$$
\begin{array}{l}
\displaystyle J_1=\int_{B_{4h}} A(\overline y-\overline
y_1)\chi_{\{t_1<t\}}F(u(y_1))\,dy_1,\\[10pt]
\displaystyle J_2=\int_{B_{4h}} A(\overline y_0-\overline
y_1)\chi_{\{t_1<t_0\}}F(u(y_1))\,dy_1,\\[10pt]
\displaystyle J_3=\int_{B_{4h}} A(\overline \eta-\overline
y_1)\chi_{\{t_1<2t_0-t\}}F(u(y_1))\,dy_1, \qquad  \overline\eta=2\overline y_0-\overline
y.
\end{array}
$$
The integrals $J_1$ and $J_2$ were already estimated in Step 2. Since $|\overline\eta-\overline y_0|=| \overline
y-\overline y_0|$, the integral $J_3$  is estimated just in the same
way as $J_1$.

The contribution of $S_h$ is estimated in the same way as in Step 2,
just using a rough estimate of the $A$'s. The contribution of $T_h$
is obviously 0.

As for  $D_h$, in this region we have, using Taylor's formula,
$$
|\mathcal{A}(y_0,y,y_1)|= |A(\overline y-\overline y_1)-2A(\overline
y_0-\overline y_1)+A(2\overline y_0-\overline y-\overline y_1)\le
\frac{ch^2}{| y_0-y_1|^4}.
$$
Since $u\in C_{\rm u}^{\alpha}(\mathbb{R}\times(0,\infty))$ for all
$\alpha\in(0,1)$, we obtain
$$
\int_{D_h} |\mathcal{A}(y_0,y,y_1)|\,|F(u(y_1))|\,dy_1 \le c h^2
\int_{D_h}\frac{dy_1}{|y_0-y_1|^{4-2\alpha}} \le ch^{2\alpha}.
$$
 In summary we get
$$
g(y_0+z)-2g(y_0)+g(y_0-z)=O(|z|^{2\alpha})
$$
for every $\alpha\in(0,1)$, $|z|<|y_0|$ (uniformly in $y_0\in\mathbb{R}\times(0,\infty)$). This estimate, together with the fact that $g$ is bounded, allows to prove  that $(-\Delta)^{\sigma/2}g(y)$ is bounded in $Q_{\tau'}$ for every  $\sigma\in(0,2)$ and $\tau'>0$. Indeed, if $\sigma\in(0,2\alpha)$ and $y\in Q_{\tau'}$, we have
$$
\begin{array}{rl}
\displaystyle|(-\Delta)^{\sigma/2}g(y)|&\displaystyle=
\left|c_\sigma\int_{\mathbb{R}^2}\frac{g(y+z)-2g(y)+g(y-z)}{|z|^{2+\sigma}}dz\right|
\\ [4mm]&\displaystyle\le
c\int_{\{|z|<\tau'\}}\frac{|z|^{2\alpha}}{|z|^{2+\sigma}}dz+
c\int_{\{|z|>\tau'\}}\frac{dz}{|z|^{2+\sigma}} \le c.
\end{array}
$$
Then, arguing in the same way as in \cite[Proposition
2.9]{Silvestre-2007} (where the boundedness of the fractional Laplacian is assumed in the whole $\mathbb{R}^2$, not only in a half-plane), if we take $\alpha\in(1/2,1)$ and $\sigma\in(1,2\alpha)$, we obtain
$g\in C^{1,\beta}(Q_{\tau'})$ for every $\beta\in(0,\sigma-1)$, with uniform
norm. We conclude that $g\in
C^{1,\alpha}_{\rm u}(Q_{\tau})$ for every $\alpha\in(0,1)$, $\tau>0$.
\qed

\subsection{$C^\infty$ regularity}  Further regularity will now be a consequence of a  result
for linear equations with smooth coefficients which has independent interest.
\begin{Theorem}
\label{th:regularity.linear.with.coefficients} Let $v$ be a bounded
weak solution to $\partial_tv+(-\Delta)^{1/2}(av+b)=0$, where the
coefficients satisfy $a,\,b\in C_{\rm u}^{1,\alpha}(\mathbb{R}\times(0,\infty))\cap
L^\infty(\mathbb{R}\times(0,\infty))$, $a(x,t)>0$. If $v\in
C_{\rm u}^\alpha(\mathbb{R}\times(0,\infty))$ then $v\in
C_{\rm u}^{1,\alpha}(\mathbb{R}\times(\tau,\infty))$ for every $\tau>0$.
\end{Theorem}

\noindent\emph{Proof. } Let $(x_0,t_0)\in\mathbb{R}^2_+$ be fixed
and denote $v_0=v(x_0,t_0)$, $a_0=a(x_0,t_0)$. Then $v$ is a
distributional solution to the inhomogeneous fractional heat
equation
$$
\partial_t v+a_0(-\Delta)^{1/2}v=(-\Delta)^{1/2}F_1+(-\Delta)^{1/2}F_2,
$$
where
$$F_1=-(a-a_0)(v-v_0), \qquad F_2=-b-v_0a.
$$
Reasoning like in the proof of Theorem~\ref{th:regularity}, we are
reduced to check that
$$
f_i(x,t)=\int_0^t\int_{\mathbb{R}}(-\Delta)^{1/2}P(x-x_1,a_0(t-t_1))F_i(x_1,t_1)\,dx_1dt_1, \qquad i=1,2
$$
are $C^{1,\alpha}$ functions, with uniform norm. It is clear that
$f_2$ inherits the  regularity of $F_2$; as to $f_1$, we use the
fact that the product $(a-a_0)(v-v_0)$ is $C^{2\alpha}$ (or
$C^{1,2\alpha-1}$ if $\alpha>1/2$) whenever $v$ is $C^{\alpha}$.
\qed

\begin{Corollary}
The strong solution to problem~\eqref{eq:main} belongs to
$C_{\rm u}^\infty(Q_\tau)$ for every $\tau>0$.
\end{Corollary}

\noindent{\it Proof. } The proof proceeds by induction. We  know that  $u\in C_{\rm u}^{1,\alpha}(Q_\tau)$, $\alpha\in(0,1)$, $\tau>0$. Assume that we have already shown that
$u\in C_{\rm u}^{k,\alpha}(Q_\tau)$ for some $k\ge1$. Then, $v_k=\partial_t^\beta\partial_x^\gamma u$, $\beta+\gamma=k$, satisfies an equation of the form $\partial_t
v_k+(-\Delta)^{1/2}(a_kv_k+b_k)=0$. Let us check that the coefficients satisfy the hypotheses of Theorem~\ref{th:regularity.linear.with.coefficients}.  On one hand, for all $k\ge1$,  $a_k=1/(1+u)$ is $C^{k,\alpha}$, hence $C^{1,\alpha}$. It is also bounded, since $u$ is nonnegative. On the other hand, as $u\in
C_{\rm u}^{k,\alpha}(Q_\tau)\cap L^\infty(Q_\tau)$, we obtain $v_k\in
C_{\rm u}^{\alpha}(Q_\tau)\cap L^\infty(Q_\tau)$. What is left is to verify that $b_k$ has the required regularity. In the case $k=1$ we have $b_1=0$, and  there is nothing to prove.  When $k=2$ we have three cases,
$$b_2=\frac{(\partial_tu)^2}{(1+u)^2},\quad \mbox{or}\quad
b_2=\frac{(\partial_xu)^2}{(1+u)^2},\quad \mbox{or}\quad b_2=\frac{\partial_tu\partial_xu}{(1+u)^2}.
$$
Since $u\in
C_{\rm u}^{2,\alpha}(Q_\tau)\cap L^\infty(Q_\tau)$, we have clearly
$b_2\in C_{\rm u}^{1,\alpha}(Q_\tau)\cap L^\infty(Q_\tau)$.
Applying Theorem~\ref{th:regularity.linear.with.coefficients}, we obtain $v_2\in C_{\rm u}^{1,\alpha}(Q_{\tau'})$, $\tau'>\tau$. Hence $u\in C_{\rm u}^{3,\alpha}(Q_{\tau'})$.

The same reasoning works for every $k\in \mathbb{N}$. Indeed, the recursion formula for the coefficients $b_k$ has the form
$$
b_{k}=\partial_ib_{k-1}+v_{k-1}\,\partial_i a\, ,
$$
where $i=x$ or $i=t$. We observe that $b_k$ is a polynomial in $\partial_t^{\beta'}\partial_x^{\gamma'}u$,
$0\le\beta'\le\beta$, $0\le\gamma'\le\gamma$, $1\le\beta'+\gamma'\le
k-1$, with coefficients involving the powers $(1+u)^{-m}$, $0<m\le k$. By the induction hypothesis, $b_{k}\in
C_{\rm u}^{1,\alpha}(Q_\tau)$. As in the step $k=2$ we conclude $u\in C_{\rm u}^{k+1,\alpha}(Q_{\tau'})$.
\qed

\subsection{Positivity}
Once the solution is regular
and the equation is satisfied in the classical sense, we can use the
Riesz representation \eqref{def-riesz} for the fractional Laplacian. Hence, at any point $(x_0,t_0)$ at which we have
$u(x_0,t_0)=0$,  we obtain
$$
\partial_tu(x_0,t_0)=\frac1\pi\mbox{ P.V.}\int_{\mathbb{R}}
\frac{\log(1+(u(s,t_0)))}{|x_0-s|^{2} }\,ds.
$$
Since $u$ is nonnegative, the right-hand side is nonnegative. Moreover, thanks to the conservation of mass, we know that  the solution $u$ is nontrivial if $f\not\equiv0$. Hence $\partial_tu(x_0,t_0)$ is strictly positive. We have thus proved the following positivity result.

\begin{Theorem}
\label{th:positivity} If $f\not\equiv0$, the solution to problem
\eqref{eq:main} is positive  for all $x\in \R$  and~$t>0$.
\end{Theorem}

\section{A nonlocal transport equation}\label{sec.trans}
\setcounter{equation}{0}

We first recall that the half-Laplacian $(-\Delta)^{1/2}$ can be
written in terms of the Hilbert transform as
$(-\Delta)^{1/2}=H\partial_x=\partial_xH$. The latter equality holds provided
that the operators are acting on a  function belonging to some $W^{1,p}(\mathbb{R})$ space,
$p>1$.

We now consider the change of variables  $(x,t,u)\mapsto (y,\tau,v)$
given by the B\"acklund type transform
$$
y=\int_0^x(1+u(s,t))\,ds-c(t),\quad \tau=t,\qquad
v(y,\tau)=\log(1+u(x,t))
$$
with $c'(t)=H(\log(1+u))(0,t)$. We denote $(y,\tau)=J(x,t)$. Notice
that the Jacobian of the transformation $J$ is
$\frac{\partial(y,\tau)}{\partial(x,t)}=1+u\ne0$, since $u\ge0$.
Then we may write the inverse
$$
x=\int_0^y e^{-v(\sigma,\tau)}\,d\sigma-\overline c(\tau),
$$
with $\overline c\,'(\tau)=-H(\log(1+u))(0,t)/(1+u(0,t))$.

We have
$$
\partial_xy=1+u,\quad \partial_ty=-H(\log(1+u))=-\widetilde H(v),
$$
where $ \widetilde H(v)=H(v\circ J)\circ J^{-1}$ is the conjugate
of the Hilbert transform $H$ by the transformation $J$.
Specifically,
\begin{equation*}\label{eq:hilbert-rare}
\begin{array}{rl}
\displaystyle\widetilde
H(v(y,\tau))=H(\log(1+u(x,t)))&\displaystyle=\frac1\pi\mbox{P.V.}
\int_{\mathbb{R}}\frac{\log(1+u(x',t))}{x-x'}\,dx'
\\ [4mm]
&\displaystyle=\frac1\pi\mbox{P.V.}
\int_{\mathbb{R}}\frac{v(y',\tau)}{\int_{y'}^ye^{v(y',\tau)-v(\sigma,\tau)}\,d\sigma}\,dy'\,.
\end{array}
\end{equation*}
With all this, equation \eqref{eq:main} becomes
\begin{equation}\label{eq:transp}
\partial_\tau v-\widetilde H(v)\,\partial_yv
+
\partial_y\widetilde H(v)=0,
\end{equation}
where $y\in\mathbb{R}$, $\tau>0$. Since we assume $u\ge0$ we get
$v\ge0$.

The $L^1$ norms of these two variables are related by
$$
\int_{\mathbb{R}}u(x,t)\,dx=\int_{\mathbb{R}}(1-e^{-v(y,\tau)})\,dy\,,
$$
and
\begin{equation*}
\int_{\mathbb{R}}v(y,\tau)\,dy=\int_{\mathbb{R}}(1+u(x,t))\log(1+u(x,t))\,dx\,.
\end{equation*}
In particular, $v_0\in L^1(\mathbb{R})$ if and only if
$u_0\in\mathcal{X}$.
On the other
hand,
$$
\partial_yv=\partial_xu, \qquad \partial_\tau v=\frac{1}{1+u}\partial_tu+ H(\log(1+u))\partial_xu.
$$
This allows to obtain regularity results for $v$ from smoothness
results for $u$. Finally, we have
$$
\begin{array}{l}
\displaystyle\int_{\mathbb{R}}|(-\Delta)^{1/4}v(y,\tau)|^2\,dy=
\int_{\mathbb{R}}|(-\Delta)^{1/4}\log(1+u(x,t))|^2\,dx, \\ [3mm]
\displaystyle\int_{\mathbb{R}}|\partial_y v(y,\tau)|^2\,dy =
\int_{\mathbb{R}}|\partial_x\log(1+u(x,t))|^2(1+u(x,t))\,dx.
\end{array}
$$

Therefore, the results of the previous sections for \eqref{eq:main}
are translated to results for~\eqref{eq:transp} as follows.

\begin{Theorem}\label{th:existence-v} Let $v_0\in
L_+^1(\mathbb{R})$. There exists a unique global in time classical
solution to equation \eqref{eq:transp} with initial value $v_0$.
\end{Theorem}

\begin{Theorem}\label{th:props-v} Let $v_0\in
L_+^1(\mathbb{R})$. The classical solution $v$ to equation
\eqref{eq:transp} with initial value $v_0$ satisfies:
\begin{enumerate}
\item[\rm (i)] \emph{$L^1$--$L^\infty$ smoothing effect:}
$\|v(\cdot,\tau)\|_\infty\le C\max\{ \tau^{-1/2}\|v_0\|_1^{1/2},\tau^{-3/4}\|v_0\|_1^{3/4}$ for all $\tau>0$;
\item[\rm (ii)] $\|v(\cdot,\tau)\|_1$ and $\|v(\cdot,\tau)\|_\infty$ are non-increasing functions of $\tau$ in $(0,\infty)$;
\item[\rm (iii)] $\displaystyle\int_{\mathbb{R}}\left(1-e^{-v(y,\tau)}\right)\,dy=
\int_{\mathbb{R}}\left(1-e^{-v_0(y)}\right)\,dy$ for every
$\tau\ge0$ \emph{(conservation law)};
\item[\rm (iv)] $v\in C^{1,\alpha}(\mathbb{R}\times(0,\infty))$ for every
$0<\alpha<1$;
\item[\rm (v)] $v(y,\tau)>0$ for every $y\in\mathbb{R}$,
$\tau>0$;
\item[\rm (vi)] $v\in L^2_{\rm loc}((0,\infty):H^1(\mathbb{R}))$.
\end{enumerate}
\end{Theorem}

\section*{Appendix A: A Nash-Trudinger inequality}
\renewcommand{\theequation}{A.\arabic{equation}}
\renewcommand{\theTheorem}{A.\arabic{Theorem}}
\renewcommand{\thesection}{A}
\label{sec.Nash-Trudinger} \setcounter{equation}{0}
\setcounter{Theorem}{0}

In this appendix we contribute a new result that falls into the
category of critical cases in embedding inequalities for spaces of
functions with weak fractional derivatives.

\begin{Theorem}\label{th:Nash-Trudinger} Let $\phi\in L^p(\mathbb{R}^N)$, $1\le p<\infty$, and assume that \
$(-\Delta)^{\gamma/2}\phi\in L^{q}(\mathbb{R}^N)$, $0<\gamma<1$,
$q=N/\gamma$. Put $r=\max\{p,\,q\}$,  $k=\lceil p/r'\rceil$, i.e.,
the least integer equal or larger than $p/r'$, $r'=r/(r-1)$.
There exists a constant $\alpha>0$ such that if $\|\phi\|_p+\|(-\Delta)^{\gamma/2}\phi\|_q\le1$ then
\begin{equation*}\label{eq:Trudinger}
\int_{\mathbb{R}^N}\left(e^{\alpha|\phi|^{r'}}-\sum_{j=0}^{k-1}\frac{(\alpha|\phi|^{r'})^j}{j!}\right)
\le 1.
\end{equation*}
\end{Theorem}

The particular case $p=q\le2$ was already proved by Strichartz in
\cite{Strichartz-1972}, using estimates on Bessel potentials. Note
that in this case  $k= 1$, so that the integrand is just
$e^{\alpha|\phi|^{p'}}-1$.  Our result covers all the possibilities
for the parameters in the critical case.

Before proceeding with the proof, we first review some related
results and preliminaries for the reader's convenience.

\noindent\textsc{Sobolev spaces of integer order. } If $1\le p<N$,
Sobolev's embedding shows that $W^{1,p}(\mathbb{R}^N)$ is
continuously embedded in $L^{r}(\mathbb{R}^N)$ for all $p\le r\le
Np/(N-p)$. If $p>N$, then $W^{1,p}(\mathbb{R}^N)$ is continuously
embedded in $L^\infty(\mathbb{R}^N)$; even more,  $\phi\in
C^{0,1-N/p}(\mathbb{R}^N)$; the same happens for $p=N=1$. The case $p=N>1$ is critical and, though
$W^{1,N}(\mathbb{R}^N)\hookrightarrow L^{r}(\mathbb{R}^N)$ for every
$1<N\le r<\infty$, it is easy to find examples of unbounded functions
in $W^{1,N}(\mathbb{R}^N)$. However, if $p=N>1$, then
\begin{equation}
\label{eq:standard.trudinger}
\int_{\mathbb{R}^N}\left(e^{\alpha|\phi|^{\frac
N{N-1}}}-\sum_{j=0}^{ N-2}\frac{\left(\alpha|\phi|^{\frac
N{N-1}}\right)^j}{j!}\right)\le 1
\end{equation}
for all $\phi$ in the unit ball of $W^{1,N}(\mathbb{R}^N)$, for some
positive $\alpha$ independent of $\phi$;  see for
example~\cite{Adams}. The proof of this result is based on the
famous analogous estimate for the case of bounded domains due to
Trudinger~\cite{Trudinger-1967}, later improved by
Moser,~\cite{Moser-1971}.

\noindent\textsc{Fractional Sobolev spaces. } Let $1\le
q<\infty$, $0<\gamma<1$.  The {\sl homogeneous
fractional Sobolev space} $\dot W^{\gamma,q}(\mathbb{R}^N)$  is defined as the completion of
$C_0^\infty(\mathbb{R})$ with the norm
$$
  \|\phi\|_{\dot{W}^{\gamma,q}}=\|(-\Delta)^{\gamma/2}\phi\|_{q}.
$$
The standard fractional Sobolev spaces are defined trough the
complete norm
$$
  \|\phi\|_{W^{\gamma,q}}=\|\phi\|_{q}+\|(-\Delta)^{\gamma/2}\phi\|_{q}.
$$
The well-known Hardy-Littlewood-Sobolev inequality
\cite{Hardy-Littlewood}, \cite{Sobolev}, states that, if
$q_*=Nq/(N-\gamma q)$, then
$$
\|\phi\|_{q_*}\le C\|(-\Delta)^{\gamma/2}\phi\|_q,
$$
for any $1<q<N/\gamma$ and $0<\gamma<1$, and thus
\begin{equation}\label{HLS1}
\dot W^{\gamma,q}(\mathbb{R}^N)\subset L^{\frac{Nq}{N-\gamma
q}}(\mathbb{R}^N).
\end{equation}

\noindent\textsc{Bessel potential spaces. } The {\sl Bessel
potential spaces}  are $L^{\gamma,q}(\mathbb{R}^N)
=\{f=J_\gamma(\phi): \phi\in L^{q}\}$, $1\le q<\infty$, where
$J_\gamma$  (the Bessel potential of order $\gamma>0$) is defined in
terms of its Fourier transform,
$$
\widehat{J_\gamma(f)}(\xi)=(1+|\xi|^2)^{-\gamma/2}\widehat{f}(\xi).
$$
These spaces were introduced by Aronszajn and
Smith~\cite{Aronszajn-Smith-1961} and
Calder\'on~\cite{Calderon-1961}, and have a natural norm,
$\|f\|_{L^{\gamma,q}}=\|J_{-\gamma}(f)\|_{L^q}.$ The space
$L^{\gamma,q}(\mathbb{R}^N)$ is equivalent to the above defined
fractional Sobolev space $W^{\gamma,q}(\mathbb{R}^N)$ for every
$0<\gamma<1$ and $1\le q<\infty$; see Stein~\cite{Stein-1961}.

\noindent\textsc{Critical Sobolev exponent. } The inclusion
\eqref{HLS1} is not valid in the critical case $q=N/\gamma$, as
pointed out in \cite{Strichartz-1967}. However, if  in addition we
know that $\phi\in L^q(\mathbb{R}^N)$, then $\phi\in
L^r(\mathbb{R}^N)$ for $N/\gamma\le r<\infty$. That is, we have the
inclusion
\begin{equation}\label{S1}
W^{\gamma,N/\gamma}(\mathbb{R}^N)\subset
L^{r}(\mathbb{R}^N),\qquad\mbox{for every } N/\gamma\le r<\infty.
\end{equation}
Indeed, in this situation $\phi$ belongs to the Bessel potential
space $L^{\gamma,N/\gamma}(\mathbb{R}^N)$, and then the result
follows from~\cite{Strichartz-1972}. Notice that the case
$r=\infty$ is not included; see \cite{Strichartz-1967}.

To go beyond the $L^r$-spaces, $q\le r<\infty$, in this critical
case $q=N/\gamma$, a careful estimate of the norms of the inclusion
\eqref{S1} using estimates of the Bessel potentials, allowed
Strichartz \cite{Strichartz-1972} to prove the inequality
\begin{equation}\label{eq:Strichartz}
\int_{\mathbb{R}^N}\left(e^{\alpha|\phi|^{\frac
N{N-\gamma}}}-1\right)\le 1
\end{equation}
for some $\alpha>0$, valid for every $\phi$ such that
$\|\phi\|_{L^{\gamma,N/\gamma}}\le1$. That is,
$L^{\gamma,N/\gamma}(\mathbb{R}^N)$ is contained in the Orlicz space
defined by the function in \eqref{eq:Strichartz}. But this result is
restricted to the range $N/2\le\gamma<N$, unless the function $\phi$
has compact support. In our case of Sobolev spaces of fractional
order $0<\gamma<1$, this means that only $N=1$ can be considered,
and then $1/2\le\gamma<1$.

\noindent\textsc{Besov spaces. } On the other hand, Peetre
\cite{Peetre} shows a restricted version of the previous inequality,
in the spirit of \eqref{eq:standard.trudinger}, valid for every
$0<\gamma<N$, for functions in the Besov space
$\Lambda^{N/\gamma,N/\gamma}_\gamma(\mathbb{R}^N)$; namely, there is  a constant $\alpha>0$ such that
\begin{equation}\label{eq:Peetre}
\int_{\mathbb{R}^N}\left(e^{\alpha|\phi|^{\frac
N{N-\gamma}}}-\sum_{j=0}^{k-1}\frac{(\alpha|\phi|^{\frac
N{N-\gamma}})^j}{j!}\right) \le 1,
\end{equation}
$k=\lceil N/(N-\gamma)\rceil$, for every $\phi$ such that
$\|\phi\|_{\Lambda^{N/\gamma,N/\gamma}_{\gamma}}\le1$.

The Besov spaces $\Lambda^{p,q}_\gamma(\mathbb{R}^N)$ are defined
through the norm
$$
  \|\phi\|_{\Lambda^{p,q}_{\gamma}}=\|\phi\|_{q}+\left(\int_{\mathbb{R}^N}\left(\int_{\mathbb{R}^N}
\frac{|\phi(x)-\phi(y)|^q}{|x-y|^{N+\gamma q}}\,dx\right)^{p/q}dy
  \right)^{1/p};
$$
see \cite{Besov-1961}.
It turns out that $\Lambda^{2,2}_\gamma(\mathbb{R}^N)=W^{\gamma,2}(\mathbb{R}^N)$. However, Besov spaces with $p=q$ and Sobolev spaces are different
whenever $q\neq2$~\cite{Stein}.

In the above-mentioned results,
the control of both the function and its derivatives,
or some quantity related to the derivatives, in the \emph{same}
$L^p$ space, yields a control in some Orlicz space. Our aim in
Theorem \ref{th:Nash-Trudinger} is  to show an  Orlicz-type estimate
analogous to \eqref{eq:standard.trudinger} and \eqref{eq:Peetre}
starting from a control of the function and its derivatives in
\emph{different} $L^p$ spaces.

We first obtain a generalization of the critical Sobolev-type
embedding~\eqref{S1},
\begin{Proposition}
Let $p\ge1$ and $0<\gamma<1$ we have
$$
L^p(\mathbb{R}^N)\cap\dot W^{\gamma,N/\gamma}(\mathbb{R}^N)\subset
L^{r}(\mathbb{R}^N)\quad\text{for every }p\le r<\infty.
$$
\end{Proposition}
\noindent{Proof. } It follows  from the Nash-Gagliardo-Nirenberg
type inequality
\begin{equation*}
   \label{eq:gagliardo.nirenberg.type.inequality}
    \|\phi\|_{rs}^r\le C(q,\gamma,N) p\|(-\Delta)^{\gamma/2}\phi\|_q
    \|\phi\|_p^{r-1},\quad r=p+1-p/q,\  s=N/(N-\gamma),
 \end{equation*}
valid for any function $\phi\in L^p(\mathbb{R}^N) \cap\dot
W^{\gamma,q}(\mathbb{R}^N)$, $p\ge1$, $q>1$,
$0<\gamma<1$, proved by the authors in~\cite{pqrv2}. Indeed, in the particular case
$q=N/\gamma$ we have $s=q'$, and thus
\begin{equation}
   \label{eq:gagliardo.nirenberg.type.inequality2}
    \|\phi\|_{p+q'}^{p+q'}\le
    Cp^{q'}\|(-\Delta)^{\gamma/2}\phi\|_q^{q'}
    \|\phi\|_p^{p}.
 \end{equation}\label{lem-NGN}
\qed


\medskip

\noindent {\sl Proof of Theorem \ref{th:Nash-Trudinger}.} As
mentioned before, the particular case $p=q\le2$ was already proved
in \cite{Strichartz-1972}. We will show how to treat
the rest of  the cases  to get a complete analysis.

\textsc{Case $p=q>2$}. It was also covered in \cite{Strichartz-1972}
under the additional restriction of asking $\phi$ to be compactly
supported. For general
functions some easy modification is needed. Indeed, for any function
$\phi\in L^{\gamma,q}(\mathbb{R}^N)$, $q=N/\gamma$, such that
$\|\phi\|_{L^{\gamma,q}}\le1$, the following estimate holds
$$
\|\phi\|_{L^r}\le A\left(1+\frac r{q'}\right)^{1/r+1/q'}
$$
for every $q\le r<\infty$, where the constant $A$ depends on $N$ and
$q$, but not on $r$, see \cite{Strichartz-1972}. We take then
$r=jq'$, $j\ge q-1$ (which implies $r\ge q$), and obtain
$$
\sum_{j\ge q-1}\frac{c^j\|\phi\|_{jq'}^{jq'}}{j!}\le \sum_{j\ge
-1}\frac{c^jA^{jq'}(j+1)^{j+1}}{j!}<\infty
$$
if we choose $c>0$ small enough. Finally, in order to have  1 in the
right-hand side of \eqref{eq:Trudinger}, we use that the function
$$
F(t)=e^{t^{q'}}-\sum_{j=0}^{k-1}\frac{t^{jq'}}{j!}
$$
satisfies $F(\lambda t)\le \lambda^{q'}F(t)$ for every $t>0$,
$0<\lambda<1$.

\textsc{Case $p< q$}. Using the Nash-Gagliardo-Nirenberg type inequality~\eqref{eq:gagliardo.nirenberg.type.inequality2} we conclude that
$\phi\in L^q(\mathbb{R}^N)$, and thus $\phi\in
L^{\gamma,q}(\mathbb{R}^N)$. We apply then the previous case.

\textsc{Case $p> q$}. The key idea is that there is a value
$0<\mu<\gamma$ such that $\phi\in L^{\mu,N/\mu}(\mathbb{R}^N)$. Indeed,
we can reach the
exponent of integration $p$ by lowering the order of differentiation. This follows from the Hardy-Littlewood-Sobolev
inclusion \eqref{HLS1}, which can be written  as
\begin{equation*}\label{HLS2}
\dot W^{\gamma_2,N/\gamma_2}(\mathbb{R}^N)\subset \dot
W^{\gamma_1,N/\gamma_1}(\mathbb{R}^N)\qquad \mbox{ for every }
0<\gamma_1<\gamma_2.
\end{equation*}
Hence, for the precise choice $\mu=N/p$ we obtain that
$(-\Delta)^{\mu/2}\phi\in L^{p}(\mathbb{R}^N)$. We may now apply the
case $p=q$ with $\gamma$ replaced by $\mu=N/p$. \qed

\noindent\emph{Remark. }   One is tempted to use the Nash-Gagliardo-Nirenberg inequality
  \eqref{eq:gagliardo.nirenberg.type.inequality2} in order to
  estimate the sum in the development of the function in
  \eqref{eq:Trudinger}. Unfortunately, the coefficient in \eqref{eq:gagliardo.nirenberg.type.inequality2}
  makes the sum divergent.

\section*{Appendix B: A calculus inequality}
\renewcommand{\theequation}{B.\arabic{equation}}
\renewcommand{\theTheorem}{B.\arabic{Theorem}}
\renewcommand{\thesection}{B}
\label{sec.Calculus} \setcounter{equation}{0}
\setcounter{Theorem}{0}

In the course of the proof of the smoothing effect we use a nice calculus inequality. Since it is not evident, we include a proof for the sake of completeness.

\begin{Lemma}
  \label{lem:exp-k} For every $x,a\ge0$ we have
  \begin{equation*}
  \label{eq:interesting.calculus.inequality}
  (e^{a x}-1)^{2}\le (e^{a}-1)(e^{a x^{2}}-1)\,.
  \end{equation*}
\end{Lemma}

\noindent{\it Proof.} We develop the function $f(x)=(e^{a}-1)(e^{a
x^{2}}-1)-(e^{a x}-1)^{2}$ in its Taylor series and rearrange the
terms as follows:
$$
\begin{array}{rl}
f(x)&\displaystyle= \sum_{n=1}^{\infty}
\frac{a^{n}}{n!}\,\sum_{k=1}^{\infty}
\frac{x^{2k}a^{k}}{k!}-\Big(\sum_{n=1}^{\infty}
\frac{x^{n}a^{n}}{n!}\Big)^{2} \\ [3mm] &\displaystyle= \sum_{n\neq
k}^{\infty} \frac{1}{n!\,k!}(a^{n+k}x^{2k}-a^{n+k}x^{n+k}) \\ [3mm]
&\displaystyle= \sum_{n\neq k}^{\infty}
\frac{1}{n!\,k!}a^{n+k}x^{n+k}(x^{k-n}-1)\,.
\end{array}
$$
By grouping the twin terms $(n,k)$ and $(k,n)$ we may restrict
ourselves to the cases $k>n$ and then
$$ f(x)\displaystyle=
\sum_{n=1}^{\infty}\sum_{k=n+1}^\infty \frac{1}{n!\,k!}
a^{n+k}x^{n+k}(x^{k-n}+x^{n-k}-2)\ge0,
$$
since the last factor is always positive for $x\ne 1$ and vanishes
for $x=1$. \qed


\


\noindent {\large\bf Acknowledgments}

\noindent FQ, AR, and JLV  partially supported by the Spanish project MTM2011-24696. AdP partially supported by
the Spanish project MTM2011-25287.

\vskip 1cm




\begin{thebibliography}{99}

\bibitem{Adams} Adams, R.~A.; Fournier, J.~J.~F.
    \lq\lq Sobolev spaces''. Second edition. Pure and Applied Mathematics (Amsterdam), 140.
    Elsevier/Academic Press, Amsterdam, 2003. ISBN: 0-12-044143-8.

\bibitem{Aronszajn-Smith-1961} Aronszajn, N.; Smith, K. T.
    \emph{Theory of Bessel potentials. I.}
    Ann. Inst. Fourier (Grenoble) 11 (1961) 385--475.

\bibitem{ac} Athanasopoulos, I.; Caffarelli, L.~A.
    \emph{Continuity of the temperature in boundary heat control problems.}
    Adv. Math. 224 (2010), no.~1, 293--315.

\bibitem{Benilan-Gariepy} B\'{e}nilan, P.; Gariepy, R.
    \emph{Strong solutions in $L^1$ of degenerate parabolic equations}.
    J. Differential Equations  119  (1995),  no.~2, 473--502.

\bibitem{Besov-1961} Besov, O.~V. \emph{Investigation of a class of function spaces in connection with
imbedding and extension theorems}. (Russian)
Trudy. Mat. Inst. Steklov. 60 (1961), 42--81.


\bibitem{cv} Caffarelli, L.~A.; Vasseur, A.
    \emph{Drift diffusion equations with fractional diffusion and the quasi-geostrophic equation.}
    Ann. of Math. (2) 171 (2010), no.~3, 1903--1930.

\bibitem{Calderon-1961} Calder\'{o}n, A.-P.
    \emph{Lebesgue spaces of differentiable functions and distributions.} 1961,
    Proc. Sympos. Pure Math., Vol. IV, pp. 33--49, American Mathematical Society, Providence, R.I.

\bibitem{CJ} Cifani, S.; Jakobsen, E.~R.
    \emph{Entropy solution theory for fractional degenerate convection-diffusion equations}.
    Ann. Inst. H. Poincar\'{e} Anal. Non Lin\'{e}aire 28 (2011), no.~3, 413--441.

\bibitem{Constantin-Lax-Majda} Constantin, P.; Lax, P.; Majda, A.
    \emph{A simple one-dimensional model for the three-dimensional vorticity.}
    Comm. Pure Appl. Math. 38 (1985), no.~6, 715--724.

\bibitem{ccf} C\'{o}rdoba, A.; C\'{o}rdoba, D.; Fontelos, M.~A.
    \emph{Formation of singularities for a transport equation with nonlocal velocity.}
    Ann. of Math. (2) 162 (2005), no.~3, 1377--1389.

\bibitem{Crandall} Crandall, M. G. \emph{An introduction to evolution governed by accretive operators}.
In \lq\lq Dynamical systems'' (Proc. Internat. Sympos., Brown Univ.,
Providence, R.I., 1974),  Vol. I, pp. 131--165. Academic Press, New
York, 1976.

\bibitem{crandall-liggett} Crandall, M.~G.; Liggett, T.~M.
\emph{Generation of semi-groups of nonlinear transformations on
general Banach spaces}.  Amer. J. Math.  93 (1971), 265--298.


\bibitem{Dong-2008} Dong, H. \emph{Well-posedness for a transport equation with
nonlocal velocity}. J. Funct. Anal. 255 (2008), no. 11, 3070--3097.

\bibitem{Hardy-Littlewood} Hardy, G. H.; Littlewood, J. E.
    \emph{Some properties of fractional integrals.~I}.
    Math. Z. 27 (1928), no.~1, 565--606.

\bibitem{Kiselev-Nazarov-Shterenberg} Kiselev, A.; Nazarov, F.; Shterenberg, R. \emph{Blow up
and regularity for fractal Burgers equation}. Dyn. Partial Differ.
Equ. 5 (2008), no.~3, 211--240.

\bibitem{Landkof} Landkof, N.~S.
    \lq\lq Foundations of modern potential theory''.
    Die Grundlehren der mathematischen Wissenschaften, Band 180. Springer-Verlag, New York-Heidelberg, 1972.

\bibitem{Moser-1971} Moser, J.
    \emph{A sharp form of an inequality by N. Trudinger}.
    Indiana Univ. Math. J. 20 (1970/71), 1077--1092.

\bibitem{OKC} Oleinik, O.~A.; Kalashnikov, A.~S.;  Czou, Y.-I.
    \emph{The Cauchy problem and boundary problems for equations of the type of non-stationary filtration}.
    Izv. Akad. Nauk SSSR. Ser. Mat. 22 (1958), 667--704. (Russian).

\bibitem{pqrv} de Pablo, A.; Quir\'os, F.; Rodr\'{\i}guez, A.; V\'azquez, J.~L.
    \emph{A fractional porous medium equation.} Adv. Math. 226 (2011), no.~2, 1378--1409.

\bibitem{pqrv2} de Pablo, A.; Quir\'os, F.; Rodr\'{\i}guez, A.; V\'azquez, J.~L.
    \emph{A general fractional porous medium equation.}
    Comm. Pure Appl. Math. 65 (2012), no.~9, 1242--1284.

\bibitem{Peetre} Peetre, J.
    \emph{Espaces d'interpolation et th\'eor\`eme de Soboleff.} (French)
    Ann. Inst. Fourier (Grenoble) 16 (1966), fasc.~1, 279--317.

\bibitem{Silvestre-2007} Silvestre, L.
    \emph{Regularity of the obstacle problem for a fractional power of the Laplace operator.}
    Comm. Pure Appl. Math. 60 (2007), no.~1, 67--112.

\bibitem{Sobolev} Sobolev, S. L.
    \emph{On a theorem of functional analysis}.
    Transl. Amer. Math. Soc. 34(2) (1963), 39--68; translation of Mat. Sb. 4 (1938) 471--497.

\bibitem{Stein-1961} Stein, E. M.
    \emph{The characterization of functions arising as potentials.}
    Bull. Amer. Math. Soc. 67 (1961) 102--104.

\bibitem{Stein} Stein, E.~M.
    \lq\lq Singular integrals and differentiability properties of functions'',
    Princeton Mathematical Series, No. 30 Princeton University Press, Princeton, N.J. 1970.

\bibitem{Strichartz-1967} Strichartz, R. S.
    \emph{Multipliers on fractional Sobolev spaces.}
    J. Math. Mech. 16 (1967) 1031--1060.

\bibitem{Strichartz-1972} Strichartz, R.~S.
    \emph{A note on Trudinger's extension of Sobolev's inequalities.}
    Indiana Univ. Math. J. 21 (1971/72), 841--842.

\bibitem{Stroock-1984} Stroock, D. W.
    \lq\lq An introduction to the theory of large deviations''.
    Universitext. Springer-Verlag, New York, 1984. ISBN: 0-387-96021-X.

\bibitem{Trudinger-1967} Trudinger, N.~S.
    \emph{On imbeddings into Orlicz spaces and some applications.}
    J. Math. Mech. 17 (1967), 473--483.

\bibitem{Varopoulos-1985} Varopoulos, N. Th.
    \emph{Hardy-Littlewood theory for semigroups.}
    J. Funct. Anal. 63 (1985), no.~2, 240--260.

\bibitem{JLVSmoothing} V\'{a}zquez, J.~L.
\lq\lq Smoothing and decay estimates for nonlinear diffusion
equations.  Equations of porous medium type''. Oxford Lecture Series
in Mathematics and its Applications, 33. Oxford University Press,
Oxford, 2006. ISBN: 978-0-19-920297-3; 0-19-920297-4.

\bibitem{vazquez} V\'{a}zquez, J.~L.
\lq\lq The porous medium equation. Mathematical theory''.    Oxford
Mathematical Monographs. The Clarendon Press, Oxford University
Press, Oxford, 2007. ISBN: 978-0-19-856903-9.


\end{thebibliography}
\end{document}